\documentclass[a4paper,11pt]{amsart}
\usepackage[margin=1in]{geometry}
\usepackage{graphicx}
\usepackage[english]{babel}
\usepackage{amsmath}
\usepackage{mathrsfs}
\usepackage{amsfonts}
\usepackage{amssymb}
\usepackage{amsthm}
\usepackage{mathtools}
\usepackage{tikz}
\usepackage{tikz-cd}
\usepackage{caption}
\usepackage{float}
\usepackage{enumerate}
\usepackage{bbm}
\usepackage{subcaption}

\usepackage[
    colorlinks=true, 
    linkcolor=blue!80!black, 
    citecolor=green!60!black
]{hyperref}
\usepackage[capitalize]{cleveref}
\usepackage{thm-restate}

\newtheorem{thm}{Theorem}[section]
\newtheorem{letterthm}{Theorem}

\newtheorem{lem}[thm]{Lemma}
\newtheorem*{lem*}{Lemma}
\newtheorem{prp}[thm]{Proposition}
\newtheorem{cor}[thm]{Corollary}

\theoremstyle{definition}
\newtheorem{dfn}[thm]{Definition}

\newtheorem{ex}[thm]{Example}
\newtheorem{rmk}[thm]{Remark}

\Crefname{lettercon}{Construction}{Constructions}
\Crefname{constructionx}{Construction}{Constructions}

\Crefname{prp}{Proposition}{Propositions}
\Crefname{cor}{Corollary}{Corollaries}
\Crefname{lettercor}{Corollary}{Corollaries}
\Crefname{letterprp}{Proposition}{Propositions}
\Crefname{lem}{Lemma}{Lemmas}
\Crefname{thm}{Theorem}{Theorems}
\Crefname{letterthm}{Theorem}{Theorems}
\Crefname{clm}{Claim}{Claims}
\Crefname{dfn}{Definition}{Definitions}
\Crefname{asn}{Assumption}{Assumptions}
\Crefname{rmk}{Remark}{Remarks}
\Crefname{ex}{Example}{Examples}
\Crefname{constr}{Construction}{Constructions}
\Crefname{assumptions}{Assumptions}{Assumptions}
\newcommand{\probref}[1]{Problem~(\ref{#1})}

\newcommand{\cB}{\mathcal{B}}
\newcommand{\cC}{\mathcal{C}}
\newcommand{\cD}{\mathcal{D}}

\newcommand{\cJ}{\mathcal{J}}
\newcommand{\cI}{\mathcal{I}}

\newcommand{\cR}{\mathcal{R}}

\newcommand{\sinf}{\mathcal{S}_\infty}
\newcommand{\cS}{\mathcal{S}}
\newcommand{\cP}{\mathcal{P}}
\newcommand{\cQ}{\mathcal{Q}}

\newcommand{\bbR}{\mathbb{R}}

\newcommand{\bbN}{\mathbb{N}}

\newcommand{\kfield}{\mathbb{K}}
\newcommand{\leqob}{\preceq_{\mathrm{ob}}}
\newcommand{\leqhom}{\preceq_{\mathrm{hom}}}
\newcommand{\cdisc}[2]{{#1_\disc#2}}
\newcommand{\rra}[2]{{#1^\rightrightarrows#2}}
\newcommand{\suff}[2]{{\mathrm{suff}(#1, #2)}}

\DeclareMathOperator{\im}{im}
\DeclareMathOperator{\coker}{coker}
\DeclareMathOperator{\Hom}{Hom}
\DeclareMathOperator{\End}{End}
\DeclareMathOperator{\Ob}{Ob}
\DeclareMathOperator{\Lan}{Lan}
\DeclareMathOperator{\length}{length}

\DeclareMathOperator*{\colim}{colim}
\let\oldlim\lim
\renewcommand{\lim}{\oldlim\nolimits}

\let\oldcolim\colim
\renewcommand{\colim}{\oldcolim\nolimits}

\newcommand{\intMod}[1]{{\Delta_{#1}\kfield}}
\newcommand{\Vect}{\mathrm{Vect}}

\DeclareMathOperator{\Mod}{Mod}

\newcommand{\id}{\mathrm{id}}
\DeclareMathOperator{\rad}{rad}

\newcommand{\defeq}{\coloneqq}  

\newcommand{\op}{\mathrm{op}}

\newcommand{\disc}{\infty}

\newcommand{\poset}[1]{\mathrm{pos}\,#1}
\newcommand{\pre}[1]{\mathrm{pre}\,#1}
\newcommand{\minCat}[2]{#1_\disc^\rightrightarrows#2}

\def\noteson{%
\gdef\samuel##1{\noindent{\color{blue}[Samuel: ##1]}}%
\gdef\justin##1{\noindent{\color{red}[Justin: ##1]}}%
\gdef\thomas##1{\noindent{\color{olive}[Thomas: ##1]}}%
\gdef\todo##1{\noindent{\color{violet}[to do: ##1]}}}

\noteson

\newcommand{\comment}[1]{}

\title[Initial Subcategories from Free Presentations]
{Initial Subcategories via Free Presentations of the Constant Functor}

\author{Thomas Br\"ustle}
\author{Justin Desrochers}
\author{Samuel Leblanc}

\AtEndDocument{\bigskip{\small%
\textsc{D\'epartement de math\'ematiques, Universit\'e de Sherbrooke\\ Sherbrooke, QC J1K 2R1, Canada} \par
\textsc{Mathematics Department, Bishop's University\\ Sherbrooke, QC J1M 1Z7, Canada} \par
  \textit{E-mail address}: 
  \texttt{tbruestl@ubishops.ca {\em and }Thomas.Brustle@USherbrooke.ca} \par
  \addvspace{\medskipamount}
  \textsc{D\'epartement de math\'ematiques, Universit\'e de Sherbrooke\\ Sherbrooke, QC J1K 2R1, Canada} \par
  \textit{E-mail address}: 
  \texttt{Justin.Desrochers@USherbrooke.ca} \par
  \addvspace{\medskipamount}
  \textsc{Department of Mathematics and Statistics, Queen's University\\ Kingston, ON K7L 3N6, Canada} \par  
  \textit{E-mail address}: 
  \texttt{Samuel.Leblanc@QueensU.ca}
}}

\hypersetup{
  pdftitle={Initial Subcategories from Free Presentations of Constant Functors},
  pdfkeywords={projective presentation, persistence module, initial functor, minimal initial subcategory, limit, Kan extension}
}

\begin{document}

\begin{abstract}
Let $\cC$ be a small category.
We characterize the full initial subcategories of $\cC$ as precisely those whose representables give a free presentation of the constant functor
$\intMod{\cC}:\cC\to\Vect$.
Such a presentation yields a formula for the limit of any functor $\cC \to \Vect$ as the kernel of a morphism between products of evaluations.

When the linearized endomorphism ring of every object is local, we show that a minimal free presentation, whenever one exists, determines the unique minimal full initial subcategory, up to equivalence.
Bass's characterization of perfect rings gives a sufficient condition for the existence of such presentations.
In particular, every finite category with no nonidentity endomorphisms, such as a finite poset, contains a unique minimal initial full subcategory.
In the poset case, our results recover the minimal initial subposet constructed in \cite{BDL26} and \cite{dey2026}. 
When $\cC$ is the free category on a finite quiver, we give an explicit construction of its minimal initial full subcategories, show that every such subcategory arises from this construction, and provide algorithms with time-complexity bounds.
\end{abstract}

\maketitle

\section{Introduction}

Limits of functors arise in many areas of mathematics, but their computation is not always efficient. 
We shall therefore be interested in functors $F:\cJ\to\cC$ that are \emph{initial}, i.e., if limits over $\cC$ can be computed by restricting, along $F$, to the category $\cJ$.
More precisely, $F$ is initial if, for all functors $M:\cC\to\cD$ 
there is a natural isomorphism
\[
\lim_{\cC} M
\;\xrightarrow{\cong}\; 
\lim_{\cJ} MF,
\]
whenever the relevant limits exist.

Of particular interest are \emph{initial subcategories}, where the initial functor $F$ is an inclusion of a subcategory $\cJ$.
If $\lim_\cC M$ is viewed as a solution to a system of equations, then $\lim_{\cJ} MF$ is the solution to a reduced system of equations.
The desire for efficient limit computation motivates two problems:

\begin{enumerate}[(1)]
    \item\label{problem1} Given a fixed small category $\cC$, 
    find an initial subcategory.
    That is, find a subcategory $\cJ$ such that the inclusion $\cJ\to \cC$ is an initial functor. 

    \item\label{problem2} Given a fixed small category $\cC$, 
    find an initial subcategory $\cI\to \cC$ which is ``minimal'' among initial subcategories of $\cC$.
\end{enumerate}

We study these problems through the category $\Mod\cC$ of modules (i.e., diagrams of vector spaces) $M:\cC\to\Vect$, for a fixed field denoted throughout by $\kfield$. 
For \probref{problem1}, we begin by considering the \emph{constant $\cC$-module} $\intMod{\cC}$ (see \cref{sec:limit-adjoint} for definition).
One way to compute the limit of $M$ in $\Mod\cC$ is via the natural isomorphism
\begin{equation*}
   \lim_{\cC}M
   \cong
   \Hom_{\Mod{\cC}}(\intMod{\cC},M).
\end{equation*}
A common technique in homological algebra is to compute the right-hand side using an exact sequence ending at $\intMod{\cC}$.
For example, we could write $\intMod{\cC}$ as the cokernel of some morphism between \emph{free modules} (coproducts of representable functors). This is a \emph{free presentation} of the constant module, and it gives a universal formula for all limits (see \cref{cor:limit-kernel}).  
The objects indexing the representable summands in such a presentation consequently determine an initial subcategory. 
\cref{thm:support-characterization} shows that these free presentations construct all full initial subcategories of $\cC$.
This relationship with free presentations was demonstrated for subposets $\cP \subseteq\bbN^d$ by Dey and Lesnick in \cite{dey2026};
they used free presentations to show the existence of initial subposets with $\Theta(m^2)$ morphisms, whenever $\cP$ is bounded below by a set of $m$ minima \cite[Theorem~3.10]{dey2026}.

For \probref{problem2}, we must first decide which notion of minimality to consider.
By \cref{thm:support-characterization}, the representable summands of a free presentation only determine the full initial subcategories.
Since full subcategories are determined by their objects, 
we seek to minimize the isomorphism classes of objects in an initial subcategory, see \cref{dfn:min-init-subcat} for a precise definition of minimality.
We discuss other notions of ``minimal'' initial subcategories in \cref{sec:related-work}.

\probref{problem2} has an especially nice solution when the endomorphism algebra of each representable $\cC$-module is \emph{local}.
In this setting, \Cref{thm:least-support} shows that a minimal free presentation of
$\intMod{\cC}$ determines a unique
minimal full initial subcategory of $\cC$.
Existence of such a presentation is automatic in several important situations considered below. 
For example, every poset satisfies the locality assumption, and the existence and uniqueness of minimal initial subcategories was shown for certain posets in \cite{BDL26} and \cite{dey2026}. 
In \cref{cor:local-existence}, we give a ring-theoretic condition for $\cC$ to have a minimal free presentation, and hence a minimal initial subcategory, using Bass's characterization of \emph{perfect rings} (see \cref{thm:bass}).

The hypothesis that representable $\cC$-modules have local endomorphism rings is very strong and often a minimal initial subcategory exists even if no minimal presentation exists. 
In the case that $\cC$ is the free category on a finite quiver, we provide a construction of initial full subcategories (see \cref{thm:construction-quiver} and \cref{sec:algorithms}) that includes all the minimal initial ones. 
If the quiver is acyclic, 
then the construction recovers the unique minimal initial subcategory obtained from a minimal presentation of the constant module.

\subsection{Organization of the Paper}
The rest of this section is dedicated to summarizing our main results and explaining the motivation for the problem coming from multiparameter persistence.
We also discuss prior work on minimal initial subposets, alternative notions of minimality, and the generalized rank invariant. 
\Cref{sec:modules} fixes the conventions for linearization, $\cC$-modules, and free presentations.
The beginning of \Cref{sec:support} recalls some preliminary results on Kan-extension and comma categories.
Then, we prove the Kan-extension and support characterization of \cref{thm:support-characterization} and
derive the universal limit formula. 
\Cref{sec:unique-min} treats minimal
presentations and explains consequences for finite posets.  
\Cref{sec:counter-examples} gives examples of the cyclic and
idempotent obstructions to unique minimal initial subcategory.  
\cref{sec:construction} gives the construction of the minimal initial full subcategory of free categories on finite quivers.
\cref{sec:examples} provides examples of this construction. 
In \cref{sec:algorithms}, we give algorithms to construct the initial full subcategories of \cref{sec:construction} and compute the time complexity. 
Finally, Appendix \ref{sec-background-rep-theory} recalls some representation theoretic facts which are used throughout the text.

We assume familiarity with category theory \cite{riehl,maclane,kashiwara05}. 
The representation theoretic background we require is explained throughout with citations. 

\begin{rmk}
Throughout, we focus on initial functors, as the results for \emph{final} functors are dual (see \cref{subsec:finality}).
In literature on final functors, such as \cite{kashiwara05} and \cite{BDL26}, we will often cite results as if they were the dual versions.
\end{rmk}

\subsection{Summary of Main Results}
\label{sec:summary}

Let $\cC$ be a small category.
Before describing our main contributions, we define the notion of minimal initial subcategory that we will consider throughout.

\begin{dfn}\label{dfn:min-init-subcat}
Define a preorder on the subcategories of a small category $\cC$ by $\cB'\leqob\cB$ if every object of $\cB'$ is isomorphic in $\cC$ to an object of $\cB$. 
A full initial subcategory $\cJ$ is \emph{minimal} if $\cJ \leqob \cI$ whenever $\cI$ is an initial subcategory such that $\cI \leqob \cJ$. 
If $\cJ\leqob \cI$ for every full initial subcategory $\cI$, then we say $\cJ$ is the \emph{unique minimal} initial full subcategory ($\cJ$ is unique up to equivalence of categories).
\end{dfn}

Note that, for full subcategories $F:\cB \hookrightarrow \cC$ and $F' : \cB' \hookrightarrow \cC$, then $\cB'\leqob\cB$ and $\cB\leqob\cB'$ if and only if there is an equivalence of categories $E : \cB \to \cB'$ such that $F \cong F'\circ E$. 
If $\cC$ is \emph{skeletal} (e.g., a poset or the free category on a quiver), then $\leqob$ is simply containment of object sets.

\begin{rmk}
Fix a skeleton $\operatorname{sk}\cC$ of a small category $\cC$.
The inclusion $\operatorname{sk}\cC \hookrightarrow \cC$ is an equivalence of categories, and therefore is an initial functor. 
Thus, one could assume without loss of generality that $\cC$ is skeletal to consider containment of objects as the partial order $\leqob$. 
\end{rmk}

\begin{rmk}
    In a category with finitely many objects, the set of full initial subcategories is finite.
    Therefore, in that case there is a minimal initial subcategory with respect to $\leqob$.
\end{rmk}

We now state our main contributions. 

There are many equivalent conditions to determine whether a functor $F:\cJ\to\cC$ is \emph{initial}, such as those in \cite[Proposition~2.5.2]{kashiwara05}.
For example, $F$ is initial if and only if the \emph{comma category} $F/ c$ is connected for every object $c$ of $\cC$.
\cref{thm:support-characterization} (\ref{item-existenceFree}) is a necessary homological condition for a functor $F:\cJ\to \cC$ to be initial.
The condition is sufficient when $F$ is fully faithful.

\begin{letterthm}\label{thm:support-characterization}
Let $F:\cJ\hookrightarrow\cC$ be the inclusion of a full subcategory. 
The following are equivalent:
\begin{enumerate}[(i)]
\item $F$ is initial.
\item The component of the counit $\Lan_F\intMod{\cJ}\to\intMod{\cC}$ is an isomorphism. 
\item $\cJ$ is the object support of some decomposition of a free presentation  of $\intMod{\cC}$.
\end{enumerate}
\end{letterthm}

We note that the equivalence of (\ref{item-initial}) and (\ref{item-canMapIso}) in \cref{thm:support-characterization} is a standard result, see for instance \cite[Thm.~4.67]{kelly}.

As an immediate consequence, for any $M\in\Mod\cC$ we obtain a formula for $\lim_\cC M$ in terms of the free presentation of $\intMod{\cC}$ in \cref{cor:limit-kernel}.

\begin{rmk}\label{rmk:support-characterization-2-implies-3}
In \cref{thm:support-characterization}, the fully faithful assumption is only used to ensure that precomposition with $F:\cJ\to \cC$ preserves representable functors. 
The proof of (\ref{item-existenceFree}) implies (\ref{item-canMapIso}) is the only claim which requires the functor to be fully faithful.
\end{rmk}

While \cref{thm:support-characterization} rephrases \probref{problem2}, we still want to find solutions. 
\cref{thm:least-support} gives sufficient conditions for the existence of 
the unique minimal full initial subcategory of $\cC$ (in the sense of \cref{dfn:min-init-subcat}).

For any objects $c,d\in\Ob\cC$, we write $\kfield\cC(c,d)$ for the vector space with basis $\cC(c,d)$.

\begin{letterthm}\label{thm:least-support}
    Suppose that $\cC$ is a small category such that $\kfield\cC(c,c)$ is a local ring for every $c\in\Ob\cC$.
    Suppose further that $\intMod{\cC}$ admits a minimal free presentation, and let $\cJ$ be any object support of such a presentation.
    Then $\cJ$ is the unique minimal initial subcategory of $\cC$.
In particular, up to equivalence of subcategories, $\cJ$ is independent of the chosen minimal free presentation and of its chosen decomposition.
\end{letterthm}

We explore consequences of this result in \cref{sec:unique-min}.
In particular, \cref{cor:local-existence} provides a criterion for the existence of a unique minimal initial subcategory.
However, the locality hypothesis in \cref{thm:least-support} is very strong: for instance, it implies that $\cJ$ contains no nonzero nonidentity idempotents.
In \cref{sec:counter-examples}, we give an example of a category with two different minimal initial subcategories for which the constant module has no minimal free presentation.

Every poset satisfies the locality assumption of \cref{thm:least-support}. 
Thus, this theorem recovers the minimal initial subposets studied in \cite{BDL26,dey2026}.
Nevertheless, locality alone does not guarantee the existence of a minimal initial subcategory. 
Indeed, \cite{BDL26} gives several examples of posets with no minimal initial subposets.

The Gabriel--Popescu embedding theorem (see \cref{prp:projectivization}) allows us to apply results and algorithms for projective presentations of modules over rings to the study of initial subcategories.
In particular, it yields the following result.

\begin{thm}
\label{thm:proj-gamma}
Let $\cC$ be a small category and let $G \defeq \bigoplus_{c\in \Ob\cC}h_c$.
Consider the functor $\Hom(G,-):\Mod\cC \to \Mod{(\operatorname{End}(G)^\op)}$.
Then every projective presentation of $\Hom(G,\intMod{\cC})$ induces a projective presentation $P_\bullet \to \intMod{\cC}$ in $\Mod\cC$.
Moreover, there exists a free presentation $F_\bullet \to \intMod{\cC}$ of which $P_\bullet$ is a direct summand.
If, in addition, $\kfield\cC(c,c)$ is local for every $c\in \Ob\cC$, then $P_\bullet$ is itself a free presentation.
\end{thm}

Algorithms for computing projective presentations are available in various settings; see \cref{rmk:gap} for some examples.
However, when $\cC$ is the free category on a finite quiver (i.e., a directed graph), we provide a direct combinatorial construction of initial full subcategories in \cref{sec:construction}.
The resulting collection of initial full subcategories includes all minimal initial full subcategories.
Our construction builds on \cite{BDL26} and applies in greater generality than \cref{thm:least-support}.

We give an informal description of the construction outlined in \cref{sec:construction}.
We define an equivalence relation on $\Ob\cC$ by mutual reachability by morphisms.
The equivalence classes form a poset $\poset{\cC}$.
Let $A\subseteq \Ob\cC$ be a complete set of representatives for the equivalence classes of $\poset{\cC}$.
We then apply the constructions of \cite{BDL26,dey2026} to obtain a minimal initial subposet of $\poset{\cC}$, which we lift to a full subcategory $\cdisc{\cC}{A}\subseteq\cC$.
Next, we adjoin the full subcategory $\rra{\cC}{A}$ consisting of objects that occur as endpoints of noncrossing paths in $\cC$.
\cref{thm:construction-quiver} shows that every initial subcategory of $\cC$ contains the objects of the full subcategory on the objects in $\minCat{\cC}{A}\defeq \cdisc{\cC}{A} \cup \rra{\cC}{A}$ for some choice of $A$. 

\begin{letterthm}\label{thm:construction-quiver}
Let $\cC$ be the free category on a finite quiver. 
\begin{enumerate}[(i)]
    \item\label{item:initialEmbed} 
    Let $A$ be a full subcategory of $\cC$ whose objects form a complete set of representatives for the equivalence classes of $\poset{\cC}$, and construct $\minCat{\cC}{A}$ as above.
    Then the fully faithful inclusion $F: \minCat{\cC}{A} \hookrightarrow \cC$ is initial.

    \item\label{item:minInitialEmbedNotBound}
    Let $G:\cJ\hookrightarrow\cC$ be an initial embedding.
    Then there exists a choice of representatives $A$ for the equivalence classes of $\poset{\cC}$ such that $\Ob\minCat{\cC}{A} \subseteq \Ob\im G$.

    \item\label{item:minInitialEmbedAcyclic}
    Suppose that $\cdisc{\cC}{A}$ is independent of the choice of representative $A$.
    If $G : \cJ \hookrightarrow \cC$ is an initial embedding, then $\Ob\minCat{\cC}{A} \subseteq \Ob\im G$ for any $A$.
    In particular, $\minCat{\cC}{A}$ is the unique minimal initial subcategory of $\cC$.
\end{enumerate}
\end{letterthm}

Suppose now that $\cC$ is the free category on a finite quiver without oriented cycles.
In this case, the construction from \cref{thm:construction-quiver} gives the unique minimal initial full subcategory and hence agrees with the subcategory computed in \cref{thm:least-support}.
We therefore obtain the following corollary as an immediate consequence of \cref{thm:construction-quiver} (\ref{item:minInitialEmbedAcyclic}) and \cref{lem:basicFactsPoset}.
\begin{cor}
Let $\cC$ be the free category on a finite quiver without oriented cycles. Then $\cdisc{\cC}{A}$ is independent of the choice of representatives $A$ for $\poset{\cC}$.
In particular, $\minCat{\cC}{A}$ is the unique minimal initial subcategory of $\cC$. 
\end{cor}

In \cref{sec:algorithms}, we provide algorithms for computing a subcategory of the form $\minCat{\cC}{A}$.
The algorithms use an arbitrary choice of representatives $A$ and do not guarantee that the resulting subcategory is minimal.
We compute $\poset{\cC}$ using standard techniques from graph theory, and compute $\cdisc{\cC}{A}$ using the algorithm of Dey and Lesnick for initial scaffolds \cite[Theorem~3.13]{dey2026}.

While \cref{thm:construction-quiver} only considers finite quivers, it has implications for infinite categories arising from \emph{thread quivers}.
Thread quivers were introduced in \cite{roosmalen14} and subsequently generalized in \cite{PRY}.
Let $Q$ be a finite quiver, and fix a collection
$\{\cP_\alpha : \alpha \text{ is an arrow in } Q\}$
of totally ordered sets.
The \emph{free category} of this thread quiver is the category $\cC$ whose morphisms are obtained by gluing together paths in the $\cP_\alpha$ at the vertices of $Q$. This construction is detailed in \cite[Section~1.2]{PRY}, where the field $k$ is replaced by the trivial ring.
It is not difficult to see that the full subcategory of $\cC$ consisting of all vertices of $Q$ is initial.
Since the composition of initial functors is initial, \cref{thm:construction-quiver} allows us to construct minimal initial subcategories of $\cC$. 
Therefore, the following is immediate.

\begin{cor}
\label{cor:thread-quivers}
Suppose that $\cC$ is the path category of a thread quiver, and let $\mathcal{V}$ be the full subcategory on the vertices of its underlying quiver.
Then every initial subcategory of $\mathcal{V}$ is also an initial subcategory of $\cC$.
In particular, \cref{thm:construction-quiver} constructs minimal initial full subcategories of $\cC$.
\end{cor}

\subsection{Persistence Motivation}
\label{sec:persistence-motivation}

Although our results are categorical, the questions we study are inspired by multiparameter persistent homology. 
A \emph{persistence module} is a functor $M:\cP \to \Vect$, where $\cP$ is a poset.
An important source of persistence modules is persistent homology: 
given a functor $T:\cP \to \mathrm{Top}$, the composite $M=H_n\circ T$ defines a persistence module, 
where $H_n:\mathrm{Top}\to \Vect$ is the $n$-th homology functor.
The functor $T$ is usually constructed from a dataset and encodes information about the shape of the data; see \cite{chazalmichel2021,botnanlesnick22,BBS26} for background.
The goal is to extract meaningful data-analytic insights from the persistence module $M$.

In multiparameter persistence, one considers posets $\cP \subseteq \bbR^d$ for $d>1$.
In this setting, the representation theory of $\cP$ is too complicated to be completely characterized by any discrete invariant \cite{carlssonzomorodian09, botnanlesnick22}.
Nevertheless, incomplete invariants still have potential for data analysis applications.

One classical invariant of a module $M : \cP \to \Vect$, where $\cP$ is connected, is the \emph{generalized rank}, defined as the rank of the canonical linear map
\[
    \lim_{\cP} M \longrightarrow \colim_{\cP} M\,.
\]
The generalized rank turns out to be equal to the \emph{multiplicity} of the constant module $\intMod{\cP}$ as a direct summand of $M$.
For $M$ taking values in finite-dimensional vector spaces, 
this equality was established by Kinser \cite{kinser2008} for connected tree posets, 
by Chambers--Letscher \cite{chambersletscher2019} for finite connected posets, 
and by Kim--M\'emoli \cite{kimmemoli2021} for locally finite posets.
The same equality was established for arbitrary functors $\cC \to \Vect$, where $\cC$ is a connected small category, in \cite{BDL26}.

Computing the generalized rank is also a key step in several other algorithms for multiparameter persistence, including GRIL \cite{GRIL} and algorithms for generalized persistence diagrams studied in \cite{kimmemoli2021,deykimmemoli2023,boo,asashibaliu2025}, among others.
Efficient computation of generalized rank motivates the following restriction problem:
\begin{enumerate}[(1)]
    \setcounter{enumi}{2}
    \item\label{problem3} Given a category $\cC$, find a smaller (or smallest) subcategory $G:\cR\hookrightarrow\cC$ such that, for every $M\in\Mod\cC$, the generalized rank of $M$ equals that of $MG$.
\end{enumerate}

One way to solve \probref{problem3} is to find a subcategory $\cR$ whose inclusion into $\cC$ is both initial and final.
Initial embeddings preserve limits, while final embeddings preserve colimits \cite[Proposition~2.5.2]{kashiwara05}; hence, an embedding with both properties preserves the generalized rank.
Conversely, if $\cP$ is a poset satisfying certain finiteness hypotheses, then any full subposet $\cR$ that solves \probref{problem3} must be both initial and final \cite[Proposition~D]{BDL26}.

Several solutions to \probref{problem3} have been developed for finitely presented modules in $\Mod\bbR^d$.
Non-full subposet solutions to \probref{problem3} have been constructed for convex connected subposets of $\bbR^2$ and $\bbR^3$ in \cite{deykimmemoli2023} and \cite{deyxin2025}, respectively, giving efficient algorithms for computing generalized rank.
For finitely presented modules in $\Mod\bbR^d$, \cite[Corollary~3.36]{dey2026} outlines an algorithm with particularly favorable asymptotic complexity.
Their algorithm finds solutions to \probref{problem2} and its dual, then combines them to obtain a solution to \probref{problem3}.

\subsection{Related Work}
\label{sec:related-work}

Our work builds primarily on two recent approaches to Problems (\ref{problem1}) through (\ref{problem3}) when the starting category is a poset. 
In \cite{BDL26}, the unique minimal initial full subposet is constructed for a broad class of posets, including finite posets and upper semilattices.
These subposets are full and are unique minimal with respect to $\leqob$. 

However, a full subcategory will usually contain more morphisms than are necessary to solve Problems (\ref{problem1}) through (\ref{problem3}).
Dey and Lesnick \cite{dey2026} therefore study \probref{problem2} under a different notion of minimality. 
They define a preorder $\leqhom$ on functors by declaring that $G\leqhom G'$, for $G:\cB\to \cC$ and $G':\cB'\to \cC$, if and only if 
\[ |\Ob \cB | \leq |\Ob \cB'| \qquad \text{and} \qquad |\Hom \cB | \leq |\Hom \cB'|.\]
They call a subcategory $\cB\subseteq\cC$ an \emph{initial scaffold} if it is minimal with respect to $\leqhom$. 
Unlike the constructions in \cite{BDL26}, initial scaffolds need not be full subcategories.
For \probref{problem2}, it suffices to consider initial subcategories, since the minima with respect to $\leqhom$ are necessarily embeddings.

Both \cite{dey2026} and \cite{BDL26} study Problems (\ref{problem1}) through (\ref{problem3}) in the context of posets, but with different approaches.
The following structure theorem for initial subposets appears in \cite{dey2026} and, under different hypotheses, in \cite{BDL26}.
We describe the construction of the object set of an initial scaffold in \cref{sec:finite-pos}. 

\begin{thm}[{\cite[Theorem~C]{BDL26}, \cite[Theorem~3.5]{dey2026}}]
    \label{thm-minimal-final-subposet}
Suppose that $\cP$ is a poset satisfying suitable finiteness conditions (for example, $\cP \subseteq \bbN^d$ is bounded below by finitely many minimal elements).
Then:
\begin{enumerate}
    \item $\cP$ admits an initial scaffold $\cS$.
    \item The inclusion $F:\cS\hookrightarrow\cP$ is initial.
    \item The image of every initial embedding $G:\cJ\to\cP$ contains an initial scaffold.
    \item\label{item:genMinimalInitial} The object set of any initial scaffold generates the unique minimal initial full subposet 
    $\cS_\infty\subseteq\cP$ with respect to $\leqob$.
\end{enumerate}
\end{thm}

Dey and Lesnick \cite{dey2026} develop efficient algorithms for constructing initial scaffolds and give detailed complexity bounds. 
Their methods apply to the classes of posets most relevant to multiparameter persistence, 
and they use initial scaffolds to obtain efficient algorithms for computing limits and generalized rank. 
The two approaches thus address closely related but distinct notions of minimality.
The constructions in \cite{BDL26} seek minimal \emph{full} subposets and give uniqueness under suitable hypotheses, 
whereas the initial scaffolds of \cite{dey2026} need not be full or unique.

The constant module $\intMod{\cP}$ is a special case of an \emph{interval module} over a poset $\cP$.
Asashiba, Gauthier, and Liu \cite{asashibaliu2025,AGL26} study the multiplicity of an interval module $I$ as a direct summand of a module $M$.
They give formulas for computing this multiplicity using free presentations of $I$, both for two-dimensional grids \cite{AGL26} and for finite posets \cite{asashibaliu2025}.
They also establish sufficient conditions on a functor $G:\cR\hookrightarrow\cP$ ensuring that the multiplicity of $IG$ in $MG$ agrees with that of $I$ in $M$.
This provides another approach to \probref{problem3} for posets.

The persistence-theoretic motivation for our work comes primarily from the literature on generalized rank, including \cite{kinser2008,chambersletscher2019,kimmemoli2021,asashibaliu2025}.
Other approaches to computing generalized rank include \cite{deykimmemoli2023,deyxin2025}, which construct poset solutions to \probref{problem3} that are minimal with respect to $\leqhom$.
These works provide the persistence-theoretic context for the restriction problems studied here; 
our focus is instead on the categorical structure.

\subsection{Acknowledgments}

The first author presented work related to \cite{BDL26} at the \emph{Maurice Auslander Centenary Lectures} \cite{thomastalk}. 
Following the presentation, C.M. Ringel asked whether projective presentations could be related to initial subposets. 
This question led us to the characterization in \cref{thm:support-characterization}.
We are very grateful to C.M. Ringel for raising this question.

TB and JD were supported by Bishop's University and Universit\'e de Sherbrooke.
TB was also supported by NSERC Discovery Grant RGPIN/2025-05047.

\subsection{Use of A.I} 
A large language model was used to review our arguments. 
Several errors were caught in this process, notably \cref{ex:counterexample-ai} was found by this A.I. 
The same A.I. also suggested the initial ideas for \cref{alg:minCat}, and for the literature review needed for \cref{alg:compute-poset-C}.

\section{Modules Over a Small Category}\label{sec:modules}

Fix a field $\kfield$ and denote by $\Vect$ the category of $\kfield$-vector spaces. 
Fix also a small category $\cC$.  
Write $\kfield\cC$ for the \emph{$\kfield$-linearization} of $\cC$. 
It has the same objects as $\cC$, while the hom-set $\kfield\cC(c,d)$ is the free vector space on the set $\cC(c,d)$ and composition is extended bilinearly. 
Put
\[
        \Mod\cC\defeq[\cC,\Vect]
\]
for the category of (covariant) functors $\cC \to \Vect$. 
We often call such functors $\cC$\emph{-modules}. 
Extending morphisms by linearity gives an equivalence
\begin{equation}\label{eq:linearization-equivalence}
 [\cC,\Vect]\cong[\kfield\cC,\Vect]_{\kfield},
\end{equation}
where the right-hand side means the category of $\kfield$-linear functors from $\kfield\cC$ to $\Vect$. 
Indeed, a functor on $\cC$ extends uniquely by linearity, and naturality on basis morphisms is equivalent to naturality on their linear combinations. 

We write $\cC(-,-)$, $\kfield\cC(-,-)$, and $\Hom_{\Mod{\cC}}(-,-)$ for the \emph{hom-bifunctors} on $\cC$, $\kfield\cC$, and $\Mod{\cC}$ respectively.
For $c\in\Ob\cC$, denote
\[
        h_c\defeq\kfield\cC(c,-).
\]
The \emph{(linear) Yoneda lemma} gives natural isomorphisms (see \cite[Section~1.9]{kelly} for details)
\begin{equation}\label{eq:linear-yoneda}
        \Hom_{\Mod{\cC}}(h_c,M)\cong Mc.
\end{equation}

The following is well-known (see, for instance \cite[Proposition~10.1.13]{prest09}), but we provide a proof for convenience. 

\begin{prp}\label{prop:free-generators}
For every $c\in\Ob\cC$, the representable module
\[
        h_c=\kfield\cC(c,-)
\]
is finitely presentable and projective. Moreover,
$\Hom_{\Mod{\cC}}(h_c,-)$ preserves arbitrary colimits.

The family $\{h_c\}_{c\in\Ob\cC}$ is a family of projective generators of
$\Mod{\cC}$.
In particular, every module is a quotient of a coproduct of representables, and the projective modules are precisely the retracts of such coproducts.
\end{prp}

\begin{proof}
The functor $\Hom_{\Mod{\cC}}(h_c,-)$
is naturally isomorphic to evaluation at $c$ by \cref{eq:linear-yoneda}. 
Since colimits in $\Mod{\cC}$ are computed pointwise,
evaluation at $c$ preserves all colimits.
In particular, it preserves filtered colimits, so $h_c$ is finitely presentable. Since evaluation also preserves epimorphisms,
$h_c$ is projective.

For any module $M$, we define an epimorphism
\begin{equation}\label{eq:canonical-free-cover}
\partial\colon
\bigoplus_{\substack{c\in\Ob\cC\\x\in M c}} h_c
 \longrightarrow M,
\end{equation}
whose summand indexed by $(c,x)$ is the natural transformation $\partial^{c,x}\colon h_c\to M$ corresponding to $x\in M c$ under the Yoneda isomorphism in \cref{eq:linear-yoneda}. 
Explicitly, for any morphism $f:c\to d$, we have $(\partial^{c,x})_d f= (M f)x$.
Taking $c=d$, $f=\id_d$, and 
$x\in M d$ shows that this morphism is pointwise surjective. 
Hence, the modules $h_c$ form a family of generators.

Finally, if $P$ is projective, applying the preceding construction to $P$ gives an epimorphism from a coproduct of representables onto $P$.
Since $P$ is projective, this epimorphism splits. 
Thus every projective is a retract of a coproduct of
representables. 
Conversely, coproducts of projectives and their retracts are projective.
\end{proof}

\begin{dfn}
We define the \emph{image} of a functor $G:\cB\to\cC$, denoted by $\im G$, as the smallest subcategory of $\cC$ containing all compositions of morphisms $G\beta$, for $\beta\in\cB(b,b')$. 
Equivalently, $\im G$ is the subcategory generated by $G(\Ob\cB)$ and $G(\Hom\cB)$.
\end{dfn}

If $G$ is injective on objects and faithful, then every morphism in $\im G$ is of the form $G\beta$ for some $\beta\in\cB(b,b')$. 
If $G$ is injective on objects and fully faithful, then $\im G$ is a full subcategory of $\cC$, and $G$ induces an isomorphism of categories $\cB\cong\im G$.

\begin{dfn}
    \label{dfn:free}
A module $P$ in $\Mod\cC$ is \emph{free} if there exist a set $A$ and a function $\sigma:A\to\Ob\cC$ such that
\[
    P\cong\bigoplus_{a\in A}h_{\sigma a}.
\]
A choice of such an isomorphism, together with the set $A$ and function $\sigma$, is called a \emph{decomposition} of the free module $P$. 
The \emph{object support} of this decomposition is the full subcategory of $\cC$ with object set $\sigma A$.

A \emph{free presentation} of $M \in \Ob{\Mod\cC}$ is an exact sequence
\begin{equation}\label{eq:free-presentation}
 P_1\xrightarrow{\partial_1}P_0\xrightarrow{\partial_0}M\longrightarrow 0
\end{equation}
where $P_0$ and $P_1$ are free modules.
We often denote such a presentation by $\partial_\bullet:P_\bullet\to M$, or simply by $P_\bullet\to M$. 
A free presentation is \emph{minimal} if $\partial_0$ and $\partial_1$ are right minimal (see \cref{sec-background-rep-theory}). 

A \emph{decomposition} of a free presentation consists of specified decompositions of $P_0$ and $P_1$.
The \emph{object support} of this decomposition is the full subcategory of $\cC$ whose object set is the union of the object supports of these two decompositions.
We say that a full subcategory $\cJ\subseteq\cC$ is \emph{an object support for $P$} if it is the object support of some decomposition of $P$.
Similarly, $\cJ$ is \emph{an object support for the free presentation $P_\bullet\to M$} if it is the object support of some decomposition of that presentation.
\end{dfn}

\begin{ex}
    Let $\cC$ be a category with a single object $c$.
    We have an isomorphism of rings $\Gamma \defeq \End_{\Mod\cC}(h_c)^\op \cong \kfield\cC(c,c) $, given by \cref{eq:linear-yoneda}.
    By \cref{prp:projectivization}, the functor
    \[\Hom_{\Mod\cC}(h_c,-):\Mod\cC\longrightarrow\Mod\Gamma\]
    is an equivalence of categories.  
    In particular, free $\cC$-modules correspond to free left $\Gamma$-modules in the usual sense.
    Under this equivalence, a decomposition of a free presentation as in \cref{eq:free-presentation}, amounts to choosing bases for the domain and codomain of $\partial_1$.
\end{ex}

\begin{rmk}
If a free module admits a decomposition into representables whose endomorphism rings are local, then its object support is unique up to equivalence of categories by \cref{lem:KRSA-free}.
However, without these locality assumptions, it is important to specify a decomposition. 
In particular, one can have $\bigoplus_{a\in A}h_{\sigma a} \cong \bigoplus_{b\in B}h_{\tau b}$, while the full subcategories on $\sigma A $ and on $\tau B$ are not equivalent.
\end{rmk}

\begin{rmk}
From \cref{prop:free-generators}, every module has a free presentation.
This follows by applying \cref{eq:canonical-free-cover} first to the module and then to its kernel. 
\end{rmk}

\section{Initial Subcategories and Object Support}\label{sec:support}

The goal of this section is to prove \cref{thm:support-characterization}.
We begin by recalling some relevant definitions and preliminary results in \cref{sec:limit-adjoint,sec:proof-thmA}.

\subsection{Limits and Constant Modules}\label{sec:limit-adjoint}
For any category $\cJ$ and vector space $V$, the \emph{constant $\cJ$-module} $\Delta_\cJ V:\cJ \to \Vect$ is the functor that sends every object to $V$ and every morphism to $\id_V$.
Since $\Vect$ has all small limits, the \emph{diagonal functor} $\Delta_\cJ: \Vect \to \Mod\cJ$ has a right adjoint $\lim_\cJ: \Mod\cJ \to \Vect$ 
whenever $\cJ$ is small (see \cite[Section~4.5]{riehl} for details).
For every $M\in\Ob{\Mod\cC}$, this adjunction gives isomorphisms
\begin{equation}\label{eq:intro-represent-limit}
\Hom_{\Mod\cC}(\intMod{\cC},M)
 \cong\Hom_{\Vect}(\kfield,\lim_{\cC}M)
 \cong\lim_{\cC}M.
\end{equation}

\subsection{Comma Categories}

Let $F:\cJ\to\cC$ be a functor and $c\in \Ob\cC$.
The \emph{comma category} $F/c$ has objects given by pairs $(j,f)$ where $f\in \cC(Fj,c)$. A morphism
$(j,f)\to(j',f')$ in $F/c$ is given by a morphism
$\alpha\in \cJ(j, j')$ such that the following commutes in $\cC$
\[
\begin{tikzcd}[sep=small]
Fj\arrow[dr,"f"'] \arrow[rr, "F\alpha"] && Fj' \arrow[dl,"f'"]\\
& c
\end{tikzcd}.
\]
As a special case, we write $\cC/c \defeq \id_\cC /c$.

\begin{ex}
Let $F:\cP \to \cQ$ be a functor between posets (i.e., an order preserving function). 
For any $q\in\Ob\cQ$, the comma category $F/q$ is isomorphic to
the full subcategory of $\cP$ on the objects
\[
    \{p\in\Ob\cP\mid Fp\leq q\}.
\]
In particular, $\cQ/q$ is the full subcategory on the objects $q'\in \Ob\cQ$ such that $q'\leq q$.
If $F$ is full, then $F/q$ is isomorphic to the full subcategory of $\cQ$ with objects
    \[  \{Fp \,|\, Fp \leq q \}\,. \]
\end{ex}

The following standard characterization of initial functors can be found, for example, in \cite[Proposition~2.5.2]{kashiwara05}.

\begin{lem}
A functor $F:\cJ\to\cC$ is initial if and only if the comma category $F/c$ is nonempty and connected for every $c\in \Ob\cC$. 
\end{lem}

\subsection{Restriction and Kan Extension}

Let $\cJ$ and $\cC$ be small categories, and let
$F:\cJ\to\cC$ be a functor.
The \emph{restriction functor}
$F^*:\Mod\cC\to\Mod\cJ$ is given by precomposition with $F$.
Since $\Vect$ admits all colimits, $F^*$ has a left adjoint, the \emph{(pointwise left) Kan extension}
$\Lan_F:\Mod\cJ \to \Mod\cC$ (for details see \cite[Corollary 6.2.6]{riehl}).

\cref{lem:lan-representable} and \cref{lem:precomp-free-is-free} give formulas for Kan extension and restriction of free modules, which will be used throughout.

\begin{lem}\label{lem:lan-representable}
For every $j\in\Ob\cJ$ there is a natural isomorphism
\[
        \Lan_F \kfield\cJ(j,-)
        \cong \kfield\cC(Fj,-).
\]
\end{lem}

\begin{proof}
For every $M\in\Mod\cC$, adjunction and \cref{eq:linear-yoneda} give natural isomorphisms
\begin{align*}
    \Hom_{\Mod\cC}\left(\Lan_F\kfield\cJ(j,-),M\right)
    &\cong
    \Hom_{\Mod\cJ}\left(\kfield\cJ(j,-),F^*M\right)\\
    &\cong MFj\\
    &\cong
    \Hom_{\Mod\cC}\left(\kfield\cC(Fj,-),M\right).
\end{align*}
The result follows from Yoneda in the functor category
(see, for instance, \cite[Proposition~2.3.1]{riehl}).
\end{proof}

\begin{rmk}
In ``classical'' representation theory, for example, in the representation theory of groups or of quivers, $\Lan_F$ is often called the \emph{induction functor} and can be expressed as a tensor product.
In that context, the above \cref{lem:lan-representable} simply means that induction sends projectives to projectives. 
\end{rmk}

\begin{lem}\label{lem:precomp-free-is-free}
Suppose $F:\cJ \to\cC$ is a fully faithful functor. 
For every $j \in \cJ$, there is an isomorphism $F^* \kfield \cC(F j, -) \cong \kfield\cJ(j,-)$. 
Moreover, $F^*$ restricts to an equivalence between the full subcategories of free modules whose object supports are contained in the image of $F$ and the full subcategory of free $\cJ$-modules.
\end{lem}

\begin{proof}
Since $F$ is fully faithful, we have 
\begin{equation*}
F^*\kfield\cC(Fj,-) = \kfield\cC(Fj,F-) \cong \kfield\cJ(j,-),
\end{equation*}
proving the first claim. 

Again, because $F$ is fully faithful, the Yoneda embedding gives 
\begin{align*}
\Hom_{\Mod\cC}\left(\kfield\cC(Fj,-),\kfield\cC(Fj',-)\right)
&\cong \kfield\cC(Fj',Fj)\\
&\cong \kfield\cJ(j',j)\\
&\cong
\Hom_{\Mod\cJ}\left(\kfield\cJ(j,-),\kfield\cJ(j',-)\right)\\
&\cong
\Hom_{\Mod\cJ}\left(F^*\kfield\cC(Fj,-),
F^*\kfield\cC(Fj',-)\right).
\end{align*}
Thus, $F^*$ is fully faithful when restricted to representable functors of the form $\Lan_F \kfield\cJ(j,-) \cong \kfield\cC(F j,-)$ and their coproducts. 
By \cite[Lemma~4.5.11]{riehl}, the components of the counit $ \Lan_F (F^* \kfield\cC(F j,-)) \to \kfield\cC(F j,-)$ are isomorphisms. 
Since both functors preserve coproducts, 
$\Lan_F$ and $F^*$ restrict to the desired equivalence.
\end{proof}

\begin{rmk}
More generally, if $F:\cJ\to\cC$ is fully faithful, then $F^*$
restricts to an equivalence between the full subcategories of free
modules whose object supports are contained in the essential image
of $F$ and the category of free $\cJ$-modules.
\end{rmk}

\subsection{Proof of \cref{thm:support-characterization}}\label{sec:proof-thmA}

The equivalence of (\ref{item-initial}) and (\ref{item-canMapIso}) in \cref{thm:support-characterization} is a standard characterization of initial functors enriched over a more general category; 
for instance \cite[Thm.~4.67]{kelly}.
The following proposition gives a direct proof for $\Vect$-enriched categories.

\begin{prp}\label{prop:constant-kan-test}
Let $F:\cJ\to\cC$ be a functor.
Then $\intMod{\cJ} = F^*\intMod{\cC}$, and the component of the counit
\begin{equation}\label{eq:constant-counit}
       \Lan_F\intMod{\cJ} = \Lan_F F^*\intMod{\cC} \to\intMod{\cC}
\end{equation}
is an isomorphism if and only if $F$ is initial.
\end{prp}

\begin{proof}
The equality $\intMod{\cJ} = F^*\intMod{\cC}$ follows directly from the definition of the constant module.

For each $c\in\Ob\cC$, there is a canonical functor
$\Pi^c:F/c\to\cJ$ sending an object $(j,Fj\to c)$ to $j\in \Ob\cJ$.
We consider the functor $\intMod{\cJ}\circ \Pi^c$ in $\Mod F/c$. 
By \cite[Theorem~6.2.1~and~Corollary~6.2.6]{riehl}, we have the leftmost isomorphism in \cref{eq:kan-colim}:
\begin{equation}
(\Lan_F\intMod{\cJ}) c \cong \colim_{F/c}(\intMod{\cJ}\circ \Pi^c) \cong \bigoplus_{\pi_0(F/c)}\kfield.
\label{eq:kan-colim}
\end{equation}
where $\pi_0(F/c)$ denotes the set of connected components of $F/c$.
The rightmost isomorphism in \cref{eq:kan-colim} follows since, in general the colimit of a $\intMod{\cI}$ is a direct sum of copies of $\kfield$, indexed by the connected components in the indexing category $\cI$ (see \cite[Corollary~2.4.4]{kashiwara05}).

Now suppose $(\Lan_F\intMod{\cJ}) c \cong (\intMod{\cC}) c = \kfield$, for all $c\in \Ob\cC$. 
Then the righthand side of \cref{eq:kan-colim} is $\kfield$, so
$\pi_0(F/c)$ is a point.
Thus, $F/c$ is connected so $F$ is initial. 

Conversely, assuming $F$ is initial, we have that $\pi_0(F/c)$ is a point.
This forces $(\Lan_F\intMod{\cJ}) c \cong \kfield \cong (\intMod{\cC}) c $,
proving the proposition.
\end{proof}

\begin{thm}[{\cref{thm:support-characterization}}]
Let $F:\cJ\hookrightarrow\cC$ be the inclusion of a full subcategory.
The following are equivalent:
\begin{enumerate}[(i)]
\item\label{item-initial} $F$ is initial.
\item\label{item-canMapIso} The component of the counit $\Lan_F\intMod{\cJ}\to\intMod{\cC}$ is an isomorphism. 
\item\label{item-existenceFree} $\cJ$ is the object support of some decomposition of a free presentation  of $\intMod{\cC}$.
\end{enumerate}
\end{thm}

\begin{proof}
The equivalence of (\ref{item-initial}) and (\ref{item-canMapIso}) is \cref{prop:constant-kan-test}.

Assume (\ref{item-canMapIso}). 
By \cref{prop:free-generators}, there is a free presentation 
\[\bigoplus_{a\in A_1} Q \to \bigoplus_{a\in A_0} Q \to \intMod{\cJ} \to 0\,,\]
where $Q = \bigoplus_{j\in \Ob\cJ} \kfield\cJ(j,-)$.
The left adjoint $\Lan_F$ commutes with colimits and, by \cref{lem:lan-representable}, sends every representable $\kfield\cJ(j,-)$ to $h_{Fj} = \kfield\cC(Fj,-)$. 
Applying $\Lan_F$ to the presentation and then using (\ref{item-canMapIso}) produces a free presentation of $\intMod{\cC}$ with object support $\cJ$.

Conversely, assume (\ref{item-existenceFree}), 
and let $\cJ$ be an object support of some free presentation $P_\bullet\to\intMod{\cC}$.
Since $F$ is fully faithful, \cref{lem:precomp-free-is-free} gives $F^*h_{Fj}\cong\kfield\cJ(j,-)$, for all representable summands $h_{Fj}$ in the decomposition of $P_\bullet$. 
Moreover, restriction preserves colimits, so applying $F^*$ gives a free presentation $F^*P_\bullet\longrightarrow F^*\intMod{\cC}=\intMod{\cJ}$.

Now apply $\Lan_F$ to this presentation.
The counit is an isomorphism on the free terms of $\Lan_F F^*P_\bullet$ by \cref{lem:lan-representable} and \cref{lem:precomp-free-is-free}.
Naturality of the counit identifies the resulting differentials with those of $\Lan_F F^* P_\bullet$. 
The induced map on cokernels is therefore an isomorphism $\Lan_F\intMod{\cJ}\to\intMod{\cC}$, proving (\ref{item-canMapIso}).
\end{proof}
 
\begin{rmk}
Like \cref{prop:constant-kan-test},
the equivalence of (\ref{item-canMapIso}) and (\ref{item-existenceFree}) can be generalized to functors enriched over a suitable monoidal category $\mathcal V$.
Indeed, the proofs of \cref{thm:support-characterization} hold if $h_c$ is replaced by the set-valued functor $\cC(c,-)$.
The enriched version of \cref{thm:support-characterization} follows by applying the \emph{free $\mathcal V$-category functor}, defined for instance in \cite[Section~2.5]{kelly}.

We choose to focus on $\Mod\cC$, as this is the most natural setting for the results in \cref{sec:unique-min}.
Moreover, it is not clear how to generalise
\cite[Theorem~3.2]{BDL26} or the limit-computation results of \cite{dey2026} to enriching categories other than $\Vect$.
\end{rmk}

\begin{cor}\label{cor:any-support-initial}
Let $P_\bullet\to\intMod{\cC}$ be a free presentation of
$\intMod{\cC}$.
Then the object support of every decomposition of this presentation
is a full initial subcategory of $\cC$.
\end{cor}

\begin{cor}\label{cor:limit-kernel}
Fix a decomposition of a free presentation
\[
 \bigoplus_{a\in A_1}h_{\sigma_1 a}\xrightarrow{\partial}
 \bigoplus_{a\in A_0}h_{\sigma_0 a}\longrightarrow
 \intMod{\cC}\longrightarrow0.
\]
For every $M\in\Mod\cC$, there is an exact sequence
\begin{equation*}
 0\longrightarrow\lim_{\cC}M\longrightarrow
 \prod_{a\in A_0}M(\sigma_0 a)\xrightarrow{\Hom(\partial,M)}
 \prod_{a\in A_1}M(\sigma_1 a).
\end{equation*}
In particular, this presentation computes $\lim_{\cC}M$ using only
the values of $M$ on objects in its object support.
\end{cor}

\begin{proof}
Apply the left-exact functor $\Hom_{\Mod\cC}(-,M)$ to the first sequence. The second sequence is obtained from this using
\cref{eq:intro-represent-limit}, \cref{eq:linear-yoneda}, and the fact that $\Hom_{\Mod\cC}(-,M)$ sends coproducts to products. 
\end{proof}

\subsection{Final Functors and Duality}\label{subsec:finality}

A functor $F : \cJ \to \cC$ is \emph{final}%
\footnote{
In some books such as \cite{kashiwara05}, final functors are instead called ``cofinal functors". We will use the terms initial and final instead.
} 
if colimits over $\cC$ can be computed by restriction along $F$. 
Initial functors are dual to final functors in the sense that
$F:\cJ \to \cC$ is final if and only if $F^\op : \cJ^\op \to \cC^\op$ is initial.

Applying \cref{thm:support-characterization} to the opposite categories gives the following dual characterization.
\begin{cor}\label{cor:dual-support}
Let $F:\cJ\hookrightarrow\cC$ be the inclusion of a full subcategory.
The following are equivalent:
\begin{enumerate}[(i)]
\item $F$ is final.
\item The component of the counit $\Lan_{F^\op}\intMod{\cJ^\op}\to\intMod{\cC^\op}$ is an isomorphism. 
\item $\cJ^\op$ is an object support for some free presentation of $\intMod{\cC^\op}$.
\end{enumerate}
\end{cor}

\begin{rmk}
Due to the previous \cref{cor:dual-support}, the correct dual notion of free presentation is again a free presentation, but on the opposite category. 
In particular, we should not replace free presentations by \emph{injective copresentations}, but instead replace $\kfield\cC(c,-)$ by $\kfield \cC^\op(c,-)$. 
\end{rmk}

\section{Categories with Unique Minimal Initial Full Subcategory}\label{sec:unique-min}

\subsection{Locality}

If $c$ is an object in a small category $\cC$ then the vector space $\kfield\cC(c,c)$ of endomorphisms is a ring with multiplication given by composition and multiplicative identity $\id_c$.
A ring is called \emph{local} if it has a unique maximal ideal.
We study consequences of \cref{thm:support-characterization} when the endomorphism rings of $\kfield\cC$ are local.

\begin{thm}[{\cref{thm:least-support}}]
    Suppose that $\cC$ is a small category such that $\kfield\cC(c,c)$ is a local ring for every $c\in\Ob\cC$.
    Suppose further that $\intMod{\cC}$ admits a minimal free presentation, and let $\cJ$ be any object support of such a presentation.
    Then $\cJ$ is the unique minimal initial subcategory of $\cC$.
In particular, up to equivalence of subcategories, $\cJ$ is independent of the chosen minimal free presentation and of its chosen decomposition.
\end{thm}

\begin{proof}
We know that $\cJ$ yields an initial full subcategory by \cref{cor:any-support-initial}. 
It remains to show that $\cJ \leqob \cJ'$, for any initial full subcategory $F:\cJ'\hookrightarrow\cC$. 
Choose a free presentation $Q_\bullet \to \intMod{\cJ'}$ and apply $\Lan_F$. 
By \cref{thm:support-characterization}, $\Lan_F Q_\bullet$ is a free presentation of $\intMod{\cC}$ with object support contained in the image of $F$.
By \cref{lem:minimal-comparison}, $P_\bullet$ is a retract of $\Lan_F Q_\bullet$. 

Fix decompositions $\Lan_F Q_n =\bigoplus_{a\in A_n}h_{F(\sigma_n a)} $, for $n=0,1$.
By \cref{prp:projective-free}, there are subsets
$A_n'\subseteq A_n$ such that
\[
P_n\cong\bigoplus_{a\in A_n'}h_{F(\sigma_n a)}.
\]
Thus $P_n$ admits a decomposition whose object support is contained
in $\cJ'$. By \cref{lem:KRSA-free}, this object support is equivalent
to the object support of the chosen decomposition of $P_n$ in the
minimal presentation. Hence every object of $\cJ$ is isomorphic in
$\cC$ to an object of $\cJ'$. Therefore
\[
\cJ\leqob\cJ'.
\]
\end{proof}

The techniques of the Gabriel--Popescu theorem (\cref{prp:projectivization}) allow tools from modules over rings to be used to study free presentations in $\Mod\cC$.

\begin{thm}[\cref{thm:proj-gamma}]
Let $\cC$ be a small category and let $G \defeq \bigoplus_{c\in \Ob\cC}h_c$.
Consider the functor $\Hom(G,-):\Mod\cC \to \Mod{(\operatorname{End}(G)^\op)}$.
Then every projective presentation of $\Hom(G,\intMod{\cC})$ induces a projective presentation $P_\bullet \to \intMod{\cC}$ in $\Mod\cC$.
Moreover, there exists a free presentation $F_\bullet \to \intMod{\cC}$ of which $P_\bullet$ is a direct summand.
If, in addition, $\kfield\cC(c,c)$ is local for every $c\in \Ob\cC$, then $P_\bullet$ is itself a free presentation.
\end{thm}

\begin{proof}
By \cref{prp:projectivization} (\ref{item:GP-usual}), $\Hom(G,-)$ has a left adjoint $L$.
By part (\ref{item:GP-weibel}) of \cref{prp:projectivization}, applying the left adjoint $L$ to a projective presentation of $\Hom(G,\intMod{\cC})$ gives a projective presentation $\partial_\bullet:P_\bullet \to \intMod{\cC}$.
Moreover, for $n=0,1$, it follows from \cref{prop:free-generators} that there exists a projectives $P_0'$ and $P_1'$ such that 
\[
P_1\oplus P_0'\oplus P_1'
\qquad\text{and}\qquad
P_0\oplus P_0'
\]
are free modules.
Then we have a free presentation
\[P_1\oplus P_0'\oplus P_1' \xrightarrow{
\begin{bmatrix}
    \partial_1 & 0&0\\
    0& \id &0
\end{bmatrix}} 
P_0\oplus P_0' \xrightarrow{\begin{bmatrix}\partial_0 & 0\end{bmatrix}} \intMod{\cC} \longrightarrow 0 \,.\]
Lastly, with the locality assumption, every projective $\cC$-module is free by \cref{prp:projective-free}. So each $P_n'$ can be taken to be $0$.
\end{proof}

From Bass's characterization of perfect rings (\cref{thm:bass}), we obtain a sufficient condition for $\cC$ to have a unique minimal initial subcategory.

\begin{cor}
\label{cor:local-existence}
Let $\cC$ be a small category and let $G \defeq \bigoplus_{c\in \Ob\cC}h_c$. 
If $\End_{\Mod\cC}(G)^\op$ is a left perfect ring then $\cC$ has finitely many objects.
Suppose, in addition, that $\kfield\cC(c,c)$ is a local ring for every $c\in\Ob\cC$.
Then $\intMod{\cC}$ has a minimal free presentation, and $\cC$ has a unique minimal initial full subcategory.
\end{cor}

\begin{proof}
Since $\Gamma\defeq \End_{\Mod\cC}(G)^\op$ is left perfect, \cref{thm:bass} implies that it has no
infinite family of nonzero pairwise orthogonal idempotents. 
The split idempotents of
\[
G=\bigoplus_{c\in\Ob\cC}h_c
\]
factoring through its summands therefore show that $\Ob\cC$ is finite.

Again by \cref{thm:bass}, every $\Gamma$-module has a projective
cover, so $\Hom(G,\intMod{\cC})$ admits a minimal projective
presentation. 
Since $\Ob\cC$ is finite,
\cref{prp:projectivization} give an equivalence between $\Mod\cC$ and $\Mod\Gamma$.
This yields a minimal projective
presentation of $\intMod{\cC}$ in $\Mod\cC$.

Furthermore, if the locality assumption is satisfied, then \cref{prp:projective-free} implies that this presentation is free. Hence $\intMod{\cC}$ admits a minimal free
presentation, and the result follows from \cref{thm:least-support}.
\end{proof}

Algorithms for computing projective presentation in special cases of \cref{cor:local-existence} are discussed in \cref{sec:algorithms}.

\subsection{Finite Posets}\label{sec:finite-pos}

\begin{dfn}\label{dfn:fin-cat}
    Call a category \emph{finite} if it is nonempty and has finitely many morphisms (hence finitely many objects). 
\end{dfn}

\begin{cor}\label{cor:finite-acyclic}
The conclusions of \cref{thm:least-support} and \cref{cor:local-existence} hold for every finite category with no nonidentity endomorphisms, and in particular for every finite poset.
\end{cor}

\begin{proof}
Since $\cC$ is finite, $\Gamma=\End(G)^{\op}$ is a
finite-dimensional $\kfield$-algebra and hence is left perfect.
Moreover,
\[
\kfield\cC(c,c)=\kfield\{\id_c\}\cong\kfield
\]
for every object $c$. Thus \cref{cor:local-existence} applies.
\end{proof}

\begin{dfn}\label{dfn:s-infty}
For a poset $\cP$, viewed as a category, define full
subposets recursively by
\begin{align*}
 \cS_0&\defeq\min\cP,\\
 \cS_n&\defeq\cS_{n-1}\sqcup
 \min\{p\in\cP:\cS_{n-1} /p
                         \text{ is not connected}\}\,,
\end{align*}
where $S_{n-1} /p$ denotes the comma category with respect to the full inclusion $S_{n-1} \hookrightarrow \cP$.
Denote $\cS_\infty \defeq \bigcup_{n\geq 0}\cS_n$. 
\end{dfn}

Assume that $\cP$ is a nonempty and finite poset. 
Because $\cP$ is finite, there exists $\ell \geq 0$ such that $\cS_\ell = \cS_{\ell+1} = \cdots = \cS_\infty$. 
The poset $\cS_\infty$ is the minimal initial full subposet as in \cref{thm-minimal-final-subposet} (\ref{item:genMinimalInitial}). 
We also have that $\Ob{\cS_\infty}$ is the object set of any initial scaffold in the sense of \cite{dey2026}, i.e.,
\begin{equation}\label{eq-deylesnick-subposet}
\cS_\infty = \{p\in \cP : (\cP/p) \setminus\{p\} \text{ is not connected}\}. 
\end{equation}

\begin{cor}[{\cite[Lemma 5.13]{dey2026}}]\label{cor:poset-recovery}
Let $\cP$ be a nonempty finite poset.
Then $\cS_\infty$ is the object support of a minimal free presentation.
\end{cor}

\begin{proof}
\cref{thm-minimal-final-subposet}, says that $\cS_\infty$ is contained in every full initial subposet and that its inclusion is initial. 
On the other hand, \cref{thm:least-support,cor:finite-acyclic} say that the support of a minimal free presentation has exactly the same minimal property. 
A poset has no distinct isomorphic objects, so the two object sets are equal.
\end{proof}

Thus, the recursive construction of $\cS_\infty$, the set in \cref{eq-deylesnick-subposet}, and the object support of a minimal free presentation are three descriptions of the same reduction. 

We shall give many examples in the next section. 

\section{Examples of Minimal Initial Subcategories}
\label{sec:counter-examples}

Recall that a \emph{congruence relation} $\cR$ on a small category $\cC$ is an equivalence relation on $\cC(c,d)$, for all $c,d\in \Ob\cC$, which is compatible with composition in $\cC$.
A congruence relation defines a \emph{quotient category} $\cC/\cR$ with the same objects as $\cC$ and morphisms given by 
\[ (\cC/\cR) (c,d) \defeq \frac{\cC(c,d)}{\cR(c,d)} \,.\]
If $\cR$ is the \emph{trivial congruence relation}, meaning two morphisms are equivalent only when they are the same, then $\cC/\cR$ is isomorphic to $\cC$. 

Let $Q$ be a (nonempty) finite quiver. 
A \emph{path} in $Q$ is a concatenation $\alpha_\ell \alpha_{\ell-1} \cdots \alpha_1$ of arrows in $Q$ such that the target of $\alpha_i$ equals the source of $\alpha_{i+1}$. 
We say that such a path has \emph{length} $\ell$. 
The \emph{free category} of $Q$ has objects the vertices of $Q$ and morphisms the set of paths. 
Composition of morphisms is concatenation of paths. 

It is clear that the free category on $Q$ is finite, in the sense of \cref{dfn:fin-cat}, if and only if $Q$ contains no cyclic paths.

\begin{ex}\label{ex:freePresPosetQuiver}
We verify the results of \cref{thm:least-support}. 
Take
\begin{equation*}
    \cP = 
    \begin{tikzcd}[sep=1em]
        6 \arrow[r] & 7 \arrow[r] & 8\\
        3 \arrow[ur, draw=none, "\circlearrowleft"{font=\normalsize, xshift=6.5pt, yshift=-6pt} description] \arrow[r] \arrow[u] & 4 \arrow[r] \arrow[u] \arrow[ur, draw=none, "\circlearrowleft"{font=\normalsize, xshift=6.5pt, yshift=-6pt} description] & 5 \arrow[u]\\
         & 1 \arrow[u]\arrow[r] \arrow[ur, draw=none, "\circlearrowleft"{font=\normalsize, xshift=6.5pt, yshift=-6pt} description] & 2\arrow[u]
    \end{tikzcd}
    \quad \text{ and }\quad 
    Q = 
    \begin{tikzcd}[sep=1em]
        6 \arrow[r] & 7 \arrow[r] & 8\\
        3 \arrow[r] \arrow[u] & 4 \arrow[r] \arrow[u] & 5 \arrow[u]\\
         & 1 \arrow[u]\arrow[r] & 2\arrow[u]
    \end{tikzcd}
\end{equation*}
a poset and a quiver, respectively. 
Let $\cC$ be the free category on $Q$.
The following are minimal free presentations in their respective categories
(this follows, for instance, from \cite[Chapter~I,~Theorem~5.8]{ASS06}):
\begin{equation*}
\text{in }\Mod \cP: 
\begin{tikzcd}
h_4 \arrow[r] & h_1 \oplus h_3 \arrow[r] & \intMod{\cP}  \arrow[r] & 0;
\end{tikzcd}
\end{equation*}
\begin{equation*}
\text{in }\Mod \cC: 
\begin{tikzcd}
h_4 \oplus h_5 \oplus h_7 \oplus h_8 \arrow[r] & h_1 \oplus h_3 \arrow[r] & \intMod{\cC} \arrow[r] & 0.
\end{tikzcd}
\end{equation*}
By \cref{thm:least-support} and \cref{cor:finite-acyclic}, we get that 
\begin{equation*}
    \begin{tikzcd}[sep=1em]
        3 \arrow[r] & 4\\
         & 1 \arrow[u]
    \end{tikzcd} \subseteq \cP
    \quad \text{ and }\quad 
    \begin{tikzcd}[sep=1em]
        & 7 \arrow[r] & 8\\
        3 \arrow[r] & 4 \arrow[r] \arrow[u] & 5 \arrow[u]\\
         & 1 \arrow[u]
    \end{tikzcd}
    \subseteq \cC
\end{equation*}
are the minimal initial full subcategories. 
Note that $\cS_\infty = \begin{tikzcd}[sep=0.7em]
        3 \arrow[r] & 4\\
         & 1 \arrow[u]
    \end{tikzcd}$, which was expected by \cref{cor:poset-recovery}. 
\end{ex}

The next few examples expose some limitations to the result of \cref{thm:least-support}.

\subsection{The Free 2-Cycle}\label{subsec:2cycle}

Let $\cC$ be the free category on
\[ Q= 
\begin{tikzcd}
1 \arrow[r, shift left, "\alpha"] & 2 \arrow[l, shift left, "\beta"] .
\end{tikzcd}
\]
The one-object full subcategories on $1$ and on $2$, 
\begin{equation*}
\cC_1 \defeq 
\begin{tikzcd}
1 \arrow[loop left, "\beta\alpha"]
\end{tikzcd}
\text{ and }
\cC_2 \defeq 
\begin{tikzcd}
2 \arrow[loop right, "\alpha\beta"]
\end{tikzcd},
\end{equation*}
are both initial. 
Indeed, any morphism from $1$ to $2$ is of the form $\alpha(\beta\alpha)^n$, and
\begin{equation*}
\begin{tikzcd}[sep=1.2em]
1 \arrow[rr, "(\beta\alpha)^{n-m}"] \arrow[dr, "\alpha(\beta\alpha)^n", swap] && 1 \arrow[dl, "\alpha(\beta\alpha)^m"]\\
&2
\end{tikzcd}
\end{equation*}
commutes whenever $m \leq n$. 
Thus, the comma category $F_1/2$ is connected, where $F_1 : \cC_1 \hookrightarrow \cC$ denotes the inclusion functor. 
It is clear that $F_1/1$ is also connected, proving that $F_1$ is initial. 
It is completely analogous for $\cC_2$. 
Both $\cC_1$ and $\cC_2$ are minimal among initial full subcategories, as they each have a single object. 

Let $0 \neq t \in \kfield$. 
For $i=1,2$, let $t_i:h_i\to\intMod{\cC}$ send every basis path to $t$. 
There are free presentations
\begin{equation}\label{eq:two-cycle-presentations}
 h_1\xrightarrow{(t-t\beta\alpha)^*}h_1
 \xrightarrow{t_1}\intMod{\cC}\longrightarrow0
 :\text{ and }:
 h_2\xrightarrow{(t-t\alpha\beta)^*}h_2
 \xrightarrow{t_2}\intMod{\cC}\longrightarrow0.
\end{equation}
Here $(t-t\beta\alpha)^*$ and $(t-t\alpha\beta)^*$ are the natural transformations corresponding, under Yoneda, to the precomposition by $t-t\beta\alpha$ and $t-t\alpha\beta$, respectively. 
These sequences are indeed exact, since using telescopic sums one shows that the kernel of $t_1$ is generated by $t-t\beta\alpha$, and similarly for $t_2$. 
\cref{thm:support-characterization} therefore confirms that $\cC_1$ and $\cC_2$ are initial subcategories. 
However, the next result tells us that no minimal free presentation exists, even if a minimal initial subcategory does exist. 

\begin{prp}\label{prop:two-cycle-no-cover}
Let $\cC$ be the free category on the 2-cycle quiver.
There is no free module $P$ that admits a right minimal epimorphism $P \twoheadrightarrow \intMod{\cC}$, and hence $\intMod{\cC}$ has no minimal free presentation.
\end{prp}

\begin{proof}
Suppose for the sake of contradiction that  $p:P\twoheadrightarrow\intMod{\cC}$ is a right minimal epimorphism, where $P$ is free and hence projective by \cref{prop:free-generators}. 
By \cref{lem:minimal-comparison} and \cref{eq:two-cycle-presentations}, there is a retraction $r_0 : h_1 \to P$ with section $s_0$, i.e., $r_0s_0 = \id_P$. 
It follows that $s_0r_0$ is an idempotent. 
But
\[
    \End(h_1)\cong\kfield[\beta\alpha]
\]
has no idempotents other than $0$ and $\id_{h_1}$. 
Since $P \neq 0$, we have $s_0r_0 \neq 0$ and therefore $P\cong h_1$. 
By Yoneda, every epimorphism $h_1\twoheadrightarrow \intMod{\cC}$ is a nonzero scalar multiple of $t_1$. 
However, the endomorphism $(\beta\alpha)^*:h_1\to h_1$ is not invertible, while $t_1(\beta\alpha)^*=t_1$. 
This contradicts right minimality. 
\end{proof}

In the next subsection, we show that the presence of a 2-cycle is not automatically a problem. 

\subsection{A 2-Cycle with a Minimal Free Presentation}\label{subsec:2cycleMinFreePres}

Add a source to the preceding quiver:
\[
\begin{tikzcd}
s \arrow[r, "\gamma"] & 1 \arrow[r, shift left, "\alpha"]
& 2 \arrow[l, shift left, "\beta"] .
\end{tikzcd}
\]
Fix again $0\neq t \in \kfield$. 
The natural transformation $t_s:h_s\to\intMod{\cC}$ has kernel generated by the differences between consecutive paths around the cycle, as shown by the next result. 

\begin{prp}\label{prop:cycle-with-cover}
The sequence
\begin{equation}\label{eq:source-cycle-presentation}
 h_1\xrightarrow{(\gamma-\beta\alpha\gamma)^*}h_s
 \xrightarrow{t_s}\intMod{\cC}\longrightarrow0
\end{equation}
is a minimal free presentation.
\end{prp}

\begin{proof}
We have
\begin{align*}
h_s(s) &\cong \kfield,\\
h_s(1) &\cong \kfield\{(\beta\alpha)^n\gamma:n\geq0\},\\
h_s(2) & \cong \kfield \{\alpha(\beta\alpha)^n\gamma : n \geq 0\}.
\end{align*}
Thus, the sequence is exact by a similar argument as in \cref{subsec:2cycle}. 
Using Yoneda (\cref{eq:linear-yoneda}), we have that $\End_{\Mod\cC}(h_s) \cong \kfield$. 
This proves that $t_s$ is right minimal. 
Because $(\gamma-\beta\alpha\gamma)^*$ is an isomorphism onto the kernel, the free presentation is minimal, thereby concluding the proof.
\end{proof}

Thus, the presence of a directed cycle does not, by itself, prevent a minimal projective presentation.
Nor is the locality assumption of \cref{thm:least-support} necessary for the existence of a minimal projective presentation, as we shall see in the next section.

\subsection{A 2-Cycle with Relations}\label{subsec:2cycle-relations}
We add relations to the quiver of \cref{subsec:2cycle}. 
Consider the quiver
\begin{equation*}
Q = 
\begin{tikzcd}
1 \arrow[r, shift left, "\alpha"] & 2 \arrow[l, shift left, "\beta"]
\end{tikzcd}\,.
\end{equation*}

Let $\cC$ be the quotient of the free category on $Q$ by the relations 
$\beta\alpha \sim \id_1$ and $\alpha\beta \sim (\alpha\beta)^2$. 
The object $1$ is initial, so $h_1\cong\intMod{\cC}$.  Consequently,
\begin{equation}
       0\longrightarrow h_1\xrightarrow{\cong}
       \intMod{\cC}\longrightarrow0
    \label{eq:2-cycle-rel-min-pres}
\end{equation}
is a minimal free presentation. 
By \cref{thm:support-characterization}, the full subcategory with object $1$ is an initial full subcategory. 
It is clearly minimal, since it only has an object. 
However, the full subcategory with object $2$ is also a minimal initial subcategory. 
Indeed, there is a single morphism from $2$ to $1$, so the comma category over $1$ is connected, 
and the comma category over $2$ is connected because every object is connected to $\id_2$. 
The category $\cC$ is finite in the sense of \cref{dfn:fin-cat}, but $\kfield\cC(2,2) = \kfield\{\id_2, e\} \cong \kfield \times \kfield$ is not local. 
Hence, we are not in the setting of \cref{thm:least-support} and \cref{cor:local-existence}. 
Moreover, the two singleton categories are not equivalent, since $\cC(1,1) = \{\id_1\}$ has a different cardinality than $\cC(2,2) = \{\id_2, e\}$. 

This example shows that a minimal free presentation need not determine a unique minimal initial full subcategory without the locality assumption of \cref{thm:least-support}.
Also, $\cC$ gives an example of a category with finite hom sets, that does not satisfy the hypothesis nor the conclusion of \cref{thm:least-support}.

\section{Constructing Minimal Initial Subquivers}
\label{sec:construction}

\begin{dfn}\label{dfn:pos-C}
Let $\cC$ be a small category.
We define an equivalence relation on $\Ob\cC$ with equivalence classes
\begin{equation*}
[c] \defeq \{c' \in \Ob\cC : \cC(c,c')\neq \emptyset \text{ and }\cC(c',c) \neq \emptyset\},
\end{equation*}
for all $c\in \Ob\cC$.
We define the \emph{poset category} of $\cC$ as the category $\poset{\cC}$ with objects $\Ob \poset{\cC} \defeq \{[c] : c \in \Ob\cC\}$
and with morphisms
\begin{equation*}
\poset{\cC}([c], [c']) \defeq 
\begin{cases}
\{[c]\to [c']\} &\text{if }\cC(c,c') \neq \emptyset,\\
\emptyset &\text{otherwise}.
\end{cases}
\end{equation*}
\end{dfn}

We collect basic facts about $\poset{\cC}$ which are obvious from the definition.

\begin{lem}\label{lem:basicFactsPoset}
Let $\cC$ be a small category.
\begin{enumerate}
    \item The category $\poset{\cC}$ is poset.
    \item Any diagram in $\cC$ yields a commutative diagram in $\poset{\cC}$. 
    \item Assume that $\cC$ is the free category on an acyclic quiver, and let $A$ be the full subcategory of $\cC$ with objects a complete set of representatives of $\Ob\poset{\cC}$. Then, $A = \cC$. 
\end{enumerate}
\end{lem}

\subsection{Extending the Construction From Posets}

We assume for the rest of this section that $\cC$ is the free category on a finite quiver. 
We define three full subcategories of $\cC$. 
Recall the definition of $\cS_\infty \subseteq \poset{\cC}$ from \cref{dfn:s-infty} and \cref{eq-deylesnick-subposet}. 
Consider a fixed complete choice of representatives of $\Ob\poset{\cC}$ and denote by $A$ the full subcategory of $\cC$ with those objects. 
Write $\cdisc{\cC}{A}$ for the full subcategory of $\cC$ with objects 
\[\Ob(\cdisc{\cC}{A})
= \{c\in\Ob A : [c]\in\cS_\infty\}.\]

Define $\rra{\cC}{A}$ to be the full subcategory of $\cC$ consisting of objects $c \in \cC \setminus \cdisc{\cC}{A}$ such that there exist $c' \in \cdisc{\cC}{A}$ and two distinct maps $f_1: c' \to c$ and $f_2 : c' \to c$ for which there is no $\overline{c} \in \cC\setminus\{c\}$ such that there is a commutative diagram
\begin{equation}\label{eq:commutingTriangles}
    \begin{tikzcd}
        c' \arrow[r, "g_1"] \arrow[dr, "f_1", swap] & \overline{c} \arrow[d, "h"] & c' \arrow[l, swap, "g_2"] \arrow[dl, "f_2"]\\
        & c
    \end{tikzcd}.
\end{equation}

We then define $\minCat{\cC}{A}$ to be the full subcategory of $\cC$ with objects 
\begin{equation*}
\Ob \minCat{\cC}{A} \defeq \Ob \cdisc{\cC}{A} \sqcup \Ob \rra{\cC}{A}.
\end{equation*} 

\begin{rmk}
The category $\cdisc{\cC}{A}$ is the analog of $\cS_\infty$ for $\cC$. 
Since all diagrams don't necessarily commute in $\cC$, as opposed to in $\poset{\cC}$, we add elements of $\rra{\cC}{A}$. 
Intuitively, $\rra{\cC}{A}$ consists of objects which are the endpoints of non-crossing disjoint paths in $\cC$. 
\end{rmk}

The rest of this section is devoted to the proof of the following result. 

\begin{thm}[{\cref{thm:construction-quiver}}] 
Let $\cC$ be the free category on a finite quiver. 
\begin{enumerate}[(i)]
    \item Let $A$ be a full subcategory of $\cC$ whose objects form a complete set of representatives for the equivalence classes of $\poset{\cC}$, and construct $\minCat{\cC}{A}$ as above.
    Then the fully faithful inclusion $F: \minCat{\cC}{A} \hookrightarrow \cC$ is initial.

    \item Let $G:\cJ\hookrightarrow\cC$ be an initial embedding.
    Then there exists a choice of representatives $A$ for the equivalence classes of $\poset{\cC}$ such that $\Ob\minCat{\cC}{A} \subseteq \Ob\im G$.

    \item Suppose that $\cdisc{\cC}{A}$ is independent of the choice of representative $A$.
    If $G : \cJ \hookrightarrow \cC$ is an initial embedding, then $\Ob\minCat{\cC}{A} \subseteq \Ob\im G$ for any $A$.
    In particular, $\minCat{\cC}{A}$ is the unique minimal initial subcategory of $\cC$.
\end{enumerate}
\end{thm}

It is clear that (\ref{item:minInitialEmbedAcyclic}) follows from (\ref{item:minInitialEmbedNotBound}). 
The proof of (\ref{item:initialEmbed}) relies on the following lemma.

\begin{lem}\label{lem:cInTwoArrows}
Let $\cC$ be the free category on a finite quiver $Q$. 
Let $c' \in \cdisc{\cC}{A}$ and $c \in \cC$. 
Suppose there are two distinct maps $f_1, f_2 \in \cC(c',c)$. 
Then, there exist $\overline{c} \in \minCat{\cC}{A}$ and a commutative diagram as in \cref{eq:commutingTriangles}.
\end{lem}

\begin{proof}
If $c \in \minCat{\cC}{A}$ then taking $f_1=g_1$, $f_2 = g_2$, and $h=\id_c$ give the desired commuting diagram. 
Suppose now that $c \notin \minCat{\cC}{A}$.
Then there must exist a $\overline{c_0} \in \cC \setminus \{c\}$ such that the diagram of \cref{eq:commutingTriangles} commutes. 
Using the notation of this diagram, we know that both $g_1$ and $g_2$ arise as distinct finite paths in $Q$. 
By definition of the free category, the morphisms in this diagram are paths in $Q$. 
Write $\length(f_n)$ for the length of $f_n$. 
Let $\ell_0 \defeq \max\{\length(f_1), \length(f_2)\}$. 
Let $\ell_1 \defeq \max\{\length(g_1), \length(g_2)\}$. 
We have $\ell_1 < \ell_0$. 
Suppose that $\overline{c_0} \notin \minCat{\cC}{A}$.  
We again get a diagram as in \cref{eq:commutingTriangles} and repeating this process yields a strictly decreasing sequence $\ell_0 > \ell_1 > \cdots$ which must terminate since $Q$ is finite.  
\end{proof}

\begin{proof}[{Proof of \cref{thm:construction-quiver} (\ref{item:initialEmbed}).}]
Let $x \in \cC$. 
Let $A$ be the full subcategory with objects given by a complete choice of representatives of $\poset{\cC}$. 
Since $F$ is a fully faithful embedding,
we identify $\im F$ with $\minCat{\cC}{A}$. 
We want to show that the comma category $F/x$ is connected. 
By \cref{thm-minimal-final-subposet}, we know that $\cS_\infty/[x]$ is connected. 
In particular, it is non-empty. 
This yields a map $f:c \to x$ for some $c$ in $\cdisc{\cC}{A} \subseteq \minCat{\cC}{A}$. 
Hence, $F/x$ is non-empty. 

Let $f : c \to x$ and $f' : c' \to x$ in $F/x$, with $c,c' \in \Ob\minCat{\cC}{A}$. 
If $x \in \minCat{\cC}{A}$, then 
\begin{equation*}
\begin{tikzcd}
c \arrow[r, "f"] \arrow[dr, "f", swap] & x \arrow[d, "\id"] & c' \arrow[l, swap, "f'"] \arrow[dl, "f'"]\\
&x
\end{tikzcd}
\end{equation*}
commutes, and so $f$ and $f'$ are connected. 
Suppose that $x \notin \minCat{\cC}{A}$. 
The fact that $ \cS_\infty /[x]$ is connected implies that there is a diagram in $\cC$:
\begin{equation}\label{eq:nonCommutativeDiag}
\begin{tikzcd}[row sep=2.4em]
c \arrow[drrr, "f", swap, bend right=15] & c_1 \arrow[drr, "f_1", swap] \arrow[r, "\gamma_1"] \arrow[l, swap, "\gamma"] & c_2 \arrow[dr, "f_2"] & \arrow[l, "\gamma_2", swap] \cdots\arrow[d, dotted] \arrow[r, "\gamma_{n-2}"] & c_{n-1}\arrow[dl, "f_{n-1}", swap] & c_n \arrow[l, "\gamma_{n-1}", swap] \arrow[dll, "f_n"] \arrow[r, "\gamma'"] & c' \arrow[dlll, "f'", bend left=15] \\
&&&x
\end{tikzcd}
\end{equation}
with $c_1,\ldots, c_n \in \cdisc{\cC}{A}$. 
This diagram need not commute in $\cC$ but its image commutes in $\poset{\cC}$.

Suppose that $f_1 \neq f\circ \gamma$.
Since $x \notin \minCat{\cC}{A}$ and using \cref{lem:cInTwoArrows}, there exists $\overline{c_1} \in \minCat{\cC}{A} = \minCat{\cC}{A}\setminus\{x\}$ with a commuting diagram, as in \cref{eq:commutingTriangles}: 
\begin{equation*}
    \begin{tikzcd}
        c_1 \arrow[r, "g_1"] \arrow[dr, "f\gamma", swap] & \overline{c_1} \arrow[d, "h"] & c_1 \arrow[l, swap, "g_2"] \arrow[dl, "f_1"]\\
        & x
    \end{tikzcd}.
\end{equation*}
We thus obtain the commutative diagram
\begin{equation*}
\begin{tikzcd}[row sep=1.3em]
c \arrow[ddrr, "f", swap, bend right=20] & c_1 \arrow[ddr, "f\gamma", swap] \arrow[l, "\gamma", swap] \arrow[r, "g_1"] & \overline{c_1} \arrow[dd, "h"] \arrow[r, equal] & \overline{c_1} \arrow[ddl, "h"] & c_1 \arrow[l, "g_2", swap] \arrow[ddll, "f_1", bend left=20]\\
&&\\
&& x
\end{tikzcd}.
\end{equation*}
Applying the same logic to $f_2, \ldots, f_n$, and $f'$, we can ``complete" the diagram of \cref{eq:nonCommutativeDiag} to obtain a commutative diagram. 
This shows that $f$ and $f'$ are connected in $F/x$, so $F$ is initial. 
\end{proof}

It remains to prove (\ref{item:minInitialEmbedNotBound}). 
We'll need the following two lemmas. 

\begin{lem}
    \label{lem:discInitialEmbed}
    Let $\cC$ be the free category on a finite quiver. 
    If $G : \cJ \hookrightarrow \cC$ is an initial embedding, then there exists a choice of representatives of $\Ob\poset{\cC}$ which yields a full subcategory $A$ of $\cC$ such that $\Ob\cdisc{\cC}{A} \subseteq \Ob\im G$.
\end{lem}

\begin{proof}
By \cref{lem:basicFactsPoset}, any diagram connecting two morphisms in $\cJ$ yields a diagram connecting these morphisms in $\poset{\cJ}$, so the induced functor $\poset{G}:\poset{\cJ} \to \poset{\cC}$ is initial. 
Denote by $\cI$ the full subcategory with $\Ob{\cI} = \Ob{\im \poset{G}}$. 
We know from \cref{thm-minimal-final-subposet} that $\cS_\infty \subseteq \cI$, which implies that there is a choice of $A$ that forces $\cdisc{\cC}{A} \subseteq \im G$. 
\end{proof}

\begin{lem}\label{lem:shortest-common-suffix}
Let $\cC$ be the free category on a quiver, let $c\in\Ob\cC$, and let
\[
    f_0,\ldots,f_n
\]
be a zigzag in the comma category $\cJ/c$. Regard the $f_\ell$ as paths
ending at $c$. If $f_m$ has minimal length among the $f_\ell$, then
$f_m$ is a suffix of every $f_\ell$.
In particular, if $f_0$ and $f_n$ have no common nonidentity suffix,
then every zigzag connecting them in $\cJ/c$ contains
$\id_c$. In particular, $c \in \cJ$.
\end{lem}

\begin{proof}
Two adjacent objects in a zigzag in $\cJ/c$ are comparable with
respect to suffixes. Indeed, a morphism
\[
   (j,f:j\to c)\longrightarrow (j',f':j'\to c)
\]
in $\cJ/c$ is a morphism $u\in \cJ(j, j')$ such that $f = f' u$, so $f'$ is a suffix of $f$.
Choose $f_m$ of minimal length. We show, moving successively to the
right and to the left in the zigzag, that $f_m$ is a suffix of every $f_\ell$.
Suppose that $f_m$ is a suffix of $f_\ell$ and that $f_\ell$ and $f_{\ell+1}$ are adjacent. 
If $f_\ell$ is a suffix of $f_{\ell+1}$, the claim follows by transitivity. 
Otherwise, $f_{\ell+1}$ is a suffix of $f_\ell$. 
Then both $f_m$ and $f_{\ell+1}$ are suffixes of the same path $f_\ell$. 
By assumption, the length of $f_m$ is at most the length of $f_{\ell+1}$, so $f_m$ is a suffix of $f_{\ell+1}$.
The same argument applies if $f_m$ is a suffix of $f_{\ell+1}$ instead, 
completing the induction.

In particular, $f_m$ is a common suffix of $f_0$ and $f_n$. 
If these two paths have no common nonidentity suffix, then necessarily $f_m=\id_c$.
\end{proof}

\begin{proof}[Proof of \cref{thm:construction-quiver} (\ref{item:minInitialEmbedNotBound})]
Again, because we only care about the objects of $\cJ$, we may assume that $G$ is fully faithful, and hence that $\im G \cong \cJ$. 
By \cref{lem:discInitialEmbed}, there is $A$ such that $\cdisc{\cC}{A} \subseteq \cJ$. 
We consider $\minCat{\cC}{A}$, the full subcategory of $\cC$ with objects $\Ob\cdisc{\cC}{A} \sqcup \Ob\rra{\cC}{A}$. 
If $\rra{\cC}{A}$ is empty, then we are done. 

Suppose otherwise and let $c \in \rra{\cC}{A}$. 
This implies that there are two distinct maps $f_1 : c' \to c$ and $f_2 : c' \to c$ in $\cC$, with $c' \in \cdisc{\cC}{A}$. 
By \cref{lem:discInitialEmbed}, $c' \in \cJ$. 
Because $G/c$ is connected, there is a zigzag in $G/c$ connecting
$f_1$ and $f_2$. 
Choose an object $h:j\to c$ of minimal path length among the objects occurring in this zigzag.
By Lemma~\ref{lem:shortest-common-suffix}, $h$ is a common suffix of
$f_1$ and $f_2$. 
Hence, there are paths $g_1,g_2:c'\to j$ such that $f_1=hg_1$ and $f_2=hg_2$. 
Thus,
\[\begin{tikzcd}
c' \arrow[r,"g_1"] \arrow[dr,"f_1"']&
j\arrow[d,"h"]&
c'\arrow[l,"g_2"']\arrow[dl,"f_2"]\\
&c&
\end{tikzcd}\]
is a commutative diagram with the same form as \cref{eq:commutingTriangles}.
Since $c\in \rra{\cC}{A}$, such a diagram exists only if $c =j $. 
Therefore, $c\in\Ob\cJ$.
\end{proof}

\subsection{Alternative Proof of \cref{thm:construction-quiver} (\ref{item:initialEmbed})}\label{sec:alt-thm-Ci} 

We conclude this section with another proof that $\minCat{\cC}{A} \hookrightarrow \cC$ is initial, which uses ideas of the previous sections. 

Let $\cC$ be the free category on a finite quiver and let $\pi:\cC \to \poset{\cC}$ be the natural projection. 
By \cref{cor:finite-acyclic}, there is a minimal free presentation 
\[ Q_1\xrightarrow{q_1}Q_0\xrightarrow{q_0}\intMod{\poset{\cC}}\longrightarrow 0. \]
By \cref{lem:quot-decomp-lift}, there exist functions $\tau_n:I_n \to \Ob\cC$ and free modules
\[ P_1 = \bigoplus_{i\in I_1} h_{\tau_1 i}, \qquad P_0 = \bigoplus_{i\in I_0} h_{\tau_0 i}\] 
such that $\Lan_\pi P_n \cong Q_n$ for $n=0,1$.

By \cref{thm:least-support}, the object support of $Q_0\oplus Q_1$ is a unique minimal initial full subcategory, and hence is equal to $\cS_\infty$, by \cref{cor:poset-recovery}.
Since $\min\poset{\cC}$ is the \emph{top} of $\intMod{\poset{\cC}}$, we have by \cite[Chapter~I, Theorem~5.8]{ASS06}, that 
\[ Q_0 \cong \bigoplus_{[a]\in\min\poset{\cC}} h_{[a]}.\] 

Fix a choice of representatives of $\poset{\cC}$ containing $\tau_1(I_1) \cup \tau_0 (I_0)$.
Let $A$ be the full subcategory of $\cC$ whose objects are these representatives, so the object support of $P_0\oplus P_1$ is $\cdisc{\cC}{A}$.
In particular, the object support of $P_0$ is $\{x\in \cdisc{\cC}{A} : [x] \in \min\poset{\cC}\}$.

We construct a free presentation with object support $\minCat{\cC}{A}$.
From \cref{lem:quot-pres-lift} and the isomorphism $\pi^*\intMod{\poset{\cC}} \cong \intMod{\cC}$, to do so it is sufficient to find a morphism $\partial_1':P_1' \to P_0$ with image
    \[ \ker (\eta_{P_0}) \cong \bigoplus_{\substack{x\in \cdisc{\cC}{A}\\ [x]\in\min\poset{\cC}}} \ker\eta_{h_{x}} \,,\]
    where $\eta:\id_{\Mod\cC} \Rightarrow \pi^* \Lan_\pi$ denotes the adjunction unit and the isomorphism follows by additivity of $\eta$ and kernels.

    For all $c\in \Ob{\cC}$ and $x \in \tau_0 (I_0)$, we fix a basis $b(x,c)$ for $(\ker \eta_{h_{x}}) c$.
    If $\ker \eta_{h_{x}} =0$ then each $b(x,c)$ is empty.
    Define a morphism
    \begin{equation}
        \psi^x: \bigoplus_{y \in \Ob{\minCat{\cC}{A}}} ~\bigoplus_{g \in b(x,y)} h_{y} ~\xlongrightarrow{} h_{x}
        \label{eq:approx-ker-minP}
    \end{equation}
    such that the restriction to the summand $h_y$ indexed by $g \in b(x,y)$ 
    is the unique morphism $g^* \,\in \Hom_{\Mod\cC}(h_{y}, h_{x})$ 
    induced by 
    the Yoneda lemma.
    Equivalently, for all $\gamma \in h_y(c)$, the morphism 
    $\psi^x$ acts via $\gamma \mapsto \gamma \circ g $.
    Now define 
    \begin{equation}
        \partial_1'\defeq \bigoplus_{i\in I_0} \psi^{\tau_0 i} 
        = \bigoplus_{\substack{x\in \cdisc{\cC}{A}\\ [x]\in\min\poset{\cC}}} \psi^{x}
        \,.\label{eq:partial-is-sum-psi}
    \end{equation}

    We show that the image of $\psi^x$ is $\ker \eta_{h_x}$, for all $x\in \tau_0(I_0)$.
    To this end, let $f_1$ and $f_2$ be distinct paths such that $f_1-f_2 \in (\ker \eta_{h_x})c$, for some $c\in \Ob{\cC}$. 
    \cref{lem:cInTwoArrows} gives a $y \in \minCat{\cC}{A}$ and a commutative diagram
\begin{equation*}
\begin{tikzcd}
x \arrow[r, "g_1"] \arrow[rd, "f_1"'] & y \arrow[d, "{h}"] & x \arrow[l, "g_2"'] \arrow[ld, "f_2"] \\
& c &
\end{tikzcd}
\end{equation*}
for some $g_1,g_2 \in \cC(x,y)$ and $h\in \cC(y, c)$.
Since $b(x,y)$ is a basis, we have $g_1-g_2 = \sum_{g\in b(x,y)} \lambda_g g$ for some $\lambda_g\in \kfield$.
By commutativity of the diagram, we have
\[ f_1-f_2 = h( g_1-g_2) = \sum_{g\in b(x,y)} \lambda_g h g= \psi^x\left(\sum_{g\in b(x,y)} \lambda_g h\right)\,,\]
where $\sum_{g\in b(x,y)} \lambda_g h$ is viewed as an element of $\bigoplus_{g \in b(x,y)} h_{y}$.
Now $(\ker \eta_{h_x})c$ is the span over $\kfield$ of $\{f_1 - f_2 :\, f_1,f_2 \in \cC(x,c) \text{ and } \pi f_1 = \pi f_2\}$. 
Thus, the image of $\psi^x$ contains $\ker \eta_{h_x}$. 
The reverse inclusion is automatic, since $\psi^x(\gamma) = \gamma \circ g$ must be in $\ker \eta_{h_x}$ if $g$ is in the kernel. 
This concludes the proof. 

\section{Examples}\label{sec:examples}

We list many examples, some of which are based on previous ones. 

\begin{ex}
Let $\cC$ be the free category on 
\begin{equation*}
Q=\begin{tikzcd}
1 \arrow[r, shift left, "\alpha"] & 2 \arrow[l, shift left, "\beta"],
\end{tikzcd}
\end{equation*}
as already discussed in \cref{subsec:2cycle}. 
We see that 
\begin{equation*}
\begin{tikzcd}
1 \arrow[loop left, "\beta\alpha"]
\end{tikzcd}
\text{ and }
\begin{tikzcd}
2 \arrow[loop left, "\alpha\beta"]
\end{tikzcd}
\end{equation*}
are the two full subcategories of $\cC$ given by the two choices of a complete set of representatives of $\poset{\cC}$. 
\end{ex}

\begin{ex}
Let $\cC$ be the free category on 
\begin{equation*}
\begin{tikzcd}
s \arrow[r, "\gamma"] & 1 \arrow[r, shift left, "\alpha"] & 2 \arrow[l, shift left, "\beta"],
\end{tikzcd}
\end{equation*}
as defined in \cref{subsec:2cycleMinFreePres}. 
As opposed to the above example, $\poset{\cC} = \{s\} \to \{1,2\}$ is not a point and $\cS_\infty = \{s\}$ yields a unique $\cC_\disc$. 
By \cref{thm:construction-quiver} (\ref{item:minInitialEmbedAcyclic}), $\cC$ admits a unique minimal initial full subcategory
\begin{equation*}
\minCat{\cC}{A} = 
\begin{tikzcd}
s \arrow[r, "\gamma"] & 1 \arrow[loop right, "\beta\alpha"]
\end{tikzcd}.
\end{equation*}
This category corresponds to the one found in \cref{subsec:2cycleMinFreePres}. 
\end{ex}

\begin{ex}
Let $\cC$ be the free category on the quiver $Q$ of  \cref{ex:freePresPosetQuiver}. 
Using the notation of that example, we have $\poset{\cC} = \cP$. 
An easy verification shows that 
\begin{equation*}
\minCat{\cC}{A} =
\begin{tikzcd}[sep=1em]
    & 7 \arrow[r] & 8\\
    3 \arrow[r] & 4 \arrow[r] \arrow[u] & 5 \arrow[u]\\
     & 1 \arrow[u]
\end{tikzcd}.
\end{equation*}
\end{ex}

\begin{ex}
The cardinality of $\minCat{\cC}{A}$ depends on the choice of $A$. 
Take $\cC$ to be the free category on 
\begin{equation*}
Q = 
\begin{tikzcd}
1 \arrow[r, shift left, bend left, "\alpha"] \arrow[r, "\beta" description] & 2 \arrow[l, shift left, bend left, "\gamma"]
\end{tikzcd}.
\end{equation*}
If $\Ob\cdisc{\cC}{A} = \{1\}$, then $\minCat{\cC}{A} = \cC$. 
However, if $\Ob\cdisc{\cC}{A} = \{2\}$, then $\minCat{\cC}{A} = \cdisc{\cC}{A}$. 
\end{ex}

\begin{ex}
\label{ex:counterexample-ai}
\cref{thm:construction-quiver} (\ref{item:initialEmbed}) is false if we consider $\cC$ to be the quotient of the free category on a quiver $Q$ by a nontrivial congruence relation $\cR$. 
Indeed, take
\begin{equation*}
Q = 
\begin{tikzcd}[sep = 1.4em]
& 1 \arrow[dl, shift right, "\alpha_1", swap] \arrow[dl, shift left, "\alpha_2"] \arrow[dr, shift left, "\beta_1"] \arrow[dr, shift right, "\beta_2", swap]\\
2 \arrow[rr, bend right, "\gamma"] && 3 \arrow[ll, bend left = 50, "\delta"]
\end{tikzcd}
\end{equation*}
and the relations are induced by $\gamma\alpha_n \sim \beta_n$ and $\delta\beta_n \sim \alpha_n$ for $n=1,2$.
Then $\sinf = \{[1]\}$, so $\cdisc{\cC}{A}$ is the full subcategory with object $1$. 
Then $2$ is not contained in $\rra{\cC}{A}$, as the paths $\alpha_n$ factor through $\delta$. Similarly $3$ is not contained in $\rra{\cC}{A}$, so $\minCat{\cC}{A} = A$.
However, the inclusion $F:A \hookrightarrow \cC$ is not initial, as $\alpha_1$ and $\alpha_2$ are not connected in $F/2$.
\end{ex}

\begin{ex}
Let $\cC$ be the quotient of the free category on a finite quiver $Q$ by a nontrivial congruence relation $\cR$. 
Even if there is a choice of $A$ which makes the statement of \cref{thm:construction-quiver} (\ref{item:initialEmbed}) hold, then it does not follow that the statement of \cref{thm:construction-quiver} (\ref{item:minInitialEmbedNotBound}) holds, even if $Q$ is acyclic. 
For example, consider the quiver
\begin{equation*}
Q = 
\begin{tikzcd}
    0\arrow[dr]\\
    1 \arrow[rr, crossing over, bend left, out=40, in=140, "\delta"]\arrow[r, "\alpha_1"] \arrow[dr, "\beta_1", swap] & 2 \arrow[d, "\gamma_1", swap, shift right] \arrow[d, "\gamma_2", shift left] \arrow[r, "\alpha_2"] & 3 \arrow[dl, "\beta_2"]\\
    & 4
\end{tikzcd}
\end{equation*}
and relations $\cR$ are induced by $\beta_1\sim \gamma_1\alpha_1$, $\gamma_2 \sim \beta_2\alpha_2$, and $\beta_1 \sim \beta_2\delta$. 
Then, 
\begin{equation*}
\poset{\cC} = 
\begin{tikzcd}[row sep=1em, column sep=1.4em]
    {[0]}\arrow[dr]\\
    {[1]} \arrow[r] & {[2]} \arrow[r] & {[3]} \arrow[r] &{[4],}
\end{tikzcd}
\end{equation*}
which implies that
\begin{equation*}
\cC_\disc = \cdisc{\cC}{A} = 
\begin{tikzcd}[row sep=1em, column sep=1.4em]
    0\arrow[dr]\\
    1 \arrow[r] & 2
\end{tikzcd}
\end{equation*}
for any choice of $A$. 
Because there are two maps $\alpha_2\alpha_1 \neq \delta: 1 \to 3$ which don't factor as in \cref{eq:commutingTriangles}, we have that $3 \in \minCat{\cC}{A}$. 
We note that we must have $3 \in \cJ$ whenever the inclusion $\cJ \hookrightarrow \cC$ is initial. 
Indeed, there is no way to connect the maps $\alpha_2\alpha_1$ and $\delta$ using only the maps connecting the objects $0$, $1$, $2$, and $4$. 
Finally, we have that $4 \in \minCat{\cC}{A}$, since the two maps $\gamma_1$ and $\gamma_2$ can't be factored as in \cref{eq:commutingTriangles}. 
However, the following diagram commutes
\begin{equation*}
\begin{tikzcd}[sep=1.6em]
2 \arrow[drr, "\gamma_1", swap] & 1 \arrow[dr, "\beta_1"] \arrow[l, "\alpha_1", swap] \arrow[rr, "\delta"] & & 3 \arrow[dl, "\beta_2", swap] & 2\arrow[dll, "\gamma_2"] \arrow[l, "\alpha_2", swap]\\
&&4
\end{tikzcd},
\end{equation*}
so $\gamma_1$ and $\gamma_2$ are connected in the comma category. 
Thus, $\minCat{\cC}{A}$ is not the minimal initial full subcategory of $\cC$. 
\end{ex}

\section{Algorithms for Finite Free Categories}\label{sec:algorithms}

In this section we give algorithms to compute $\cdisc{\cC}{A}$ and $\minCat{\cC}{A}$, when $\cC$ is the free category of a finite quiver.
Dey and Lesnick \cite{dey2026} give an algorithm which, given the Hasse quiver $(V,E)$ of a finite poset $\cP$, computes a minimal initial scaffold in time $O(|V||E|)$, using depth-first search.
If $\cC$ is a category such that $\poset{\cC}$ (see \cref{sec:construction} for the definition) is finite, then applying their algorithm to $\poset{\cC}$ gives an algorithm for computing $\cdisc{\cC}{A}$.
We combine their algorithm with existing algorithms for strongly connected components and transitive reduction to obtain the following.

Before giving the algorithms, we define several notions from the graph algorithm literature. 
The \emph{reachability preorder} $\pre\cC$ is defined by $c\preceq c'$ if and only if $\cC(c,c')\neq\emptyset$.
We have an associated equivalence relation defined by $[c]= [c']$ if and only if $c\preceq c'$ and $c'\preceq c$. 
This is the equivalence relation that defines $\poset{\cC}$. 
We call the equivalence classes (the elements of $\poset{\cC}$) the \emph{strongly connected components} of $\cC$. 
The quiver obtained from $Q$ by collapsing each strongly connected component to a single vertex is the \emph{condensation} of $Q$. The condensation is always acyclic and never contains parallel arrows. 
For a finite acyclic quiver $Q$, its \emph{transitive reduction} is the quiver with the same vertices, but with an arrow $x \to y$ if and only if there is no directed path from $x$ to $y$ of length greater than one.
Thus, the transitive reduction of the condensation of $Q$ gives the Hasse quiver of the poset $\poset{\cC}$.

\begin{prp}\label{alg:compute-poset-C}
    Let $\cC$ be the free category on a finite quiver $Q$ with $n$ vertices and $m$ arrows.
The Hasse quiver of $\poset{\cC}$ can be computed in time $O(n + m p)$, where $p = |\!\Ob(\poset{\cC})|$.
\end{prp}

\begin{proof}
By the above discussion, it suffices to compute the transitive reduction of the condensation of $Q$. 
The strongly connected components of the quiver $Q$  can be computed in time $O(m+n)$, for example using Tarjan's algorithm \cite{tarjan}. The condensation of $Q$ can therefore be constructed in time $O(m+n)$.
The condensation has $p$ vertices and at most $m$ arrows. Since it is acyclic, its transitive reduction can be computed in time $O(m p)$ using the algorithm of Goral\v{c}\'{\i}kov\'{a} and Koubek \cite{transReduct}. 
Therefore, the total running time is
\[O(m+n)+O(m+n)+O(m p)=O(n+m+m p)=O(n+m p).\]
\end{proof}

\begin{prp}\label{alg:minCat}
    Let $\cC$ be the free category on a finite quiver with $n$ vertices and $m$ arrows. 
    There is an algorithm which computes some $\minCat{\cC}{A}$ in time $O((n+m)p)$, where $p=|\!\Ob\poset{\cC}|$.
\end{prp}

\begin{proof}
    
    We claim that 
    $c\in\rra{\cC}{A}$ if and only if $c \in \Ob\cC \setminus \Ob\cdisc{\cC}{A}$ and $2 \leq \max \{ |\suff{c'}{c}| : c' \in \cdisc{\cC}{A} \} $, 
    where
    \[ \suff{c'}{c} \defeq \{ 
        \text{arrows } \alpha: x \to c \text{ in } Q 
        ~|~ [c'] \leq [x]
    \}\]
    is the set of suffixes of paths in $\cC(c',c)$.
    Indeed, 
    by \cref{lem:shortest-common-suffix}, 
    a diagram of the form in \cref{eq:commutingTriangles} exists if and only if $f_1,f_2$ have a common suffix $h$.
    Thus, $ \suff{c'}{c}$ contains distinct arrows $\alpha_1:x_1\to c$ and $\alpha_2:x_2\to c$ if and only if there are paths 
    $c' \xrightarrow{g_1} x_1 \xrightarrow{\alpha_1} c$ and $c' \xrightarrow{g_2} x_2 \xrightarrow{\alpha_2} c$ without a common suffix.

    We first compute and store the strongly connected components of the quiver $Q$.
    This is used to compute $\poset{\cC}$ anyway.
    We then complete the computation of the Hasse quiver $(V,E)$ of $\poset{\cC}$ using \cref{alg:compute-poset-C}. 
    This requires $O(n+m |V|)$.

    Dey and Lesnick \cite[Theorem~3.13]{dey2026} give an algorithm which, given the Hasse quiver $(V,E)$ of a finite poset $\cP$, computes a minimal initial scaffold in time $O(|V||E|)$.
    The object set of this initial scaffold is $\cS_\infty$, so we obtain $\Ob\cdisc{\cC}{A}$ by choosing any set of representatives $A$ for the strongly connected components.

We now describe the computation of $\rra{\cC}{A}$:
\begin{enumerate}
\item
Create an array $\texttt{CurrentMax}$ indexed by $c\in \Ob\cC$.
We initialize the value $\texttt{CurrentMax}(c)$ of $\texttt{CurrentMax}$ at each index $c$, to zero.

\item
For each $d\in\Ob\cdisc{\cC}{A}$, compute all the values $\suff{d}{c}$ as
follows.
\begin{enumerate}
    \item 
    Create an array $\texttt{NumSuff}$ indexed by $c\in \Ob\cC$, (re)initialized to zero at each index.
    This takes $O(n)$ time.

    \item 
    Create an array $\texttt{GreaterThan}_d$ indexed by $\Ob\cC$, (re)initialized to zero at each index.
    Perform an ordinary depth first search of $\poset{\cC}$, starting at $[d]$.
    At each $[d] \leq [x]$, we set the value of $\texttt{GreaterThan}_d$ to \texttt{true} for all vertices in the strongly connected component $[x]$.
    Afterwards, the value of $\texttt{GreaterThan}_d$ at index $x\in \Ob\cC$ is \texttt{true} if $[d] \leq [x]$ and \texttt{false} otherwise.
    This can be done in time $O(n+|V|+|E|)$.
    
    \item 
    For each arrow $\alpha:x \to c$ in $Q$, we increment $\texttt{NumSuff}(c)$ if $[d] \leq [x]$; i.e., if the value of $\texttt{GreaterThan}_d$ at index $c$ is $\texttt{true}$.
    Checking the entries $\texttt{GreaterThan}_d$ requires constant time  
    and the cost to search all arrows is $O(m)$. 
    
    \item 
    After iterating through all arrows, we have $\texttt{NumSuff}(c) = |\suff{d}{c}|$, for all $c\in \Ob\cC$. 
    If $\texttt{CurrentMax}(c) \leq \texttt{NumSuff}(c)$ then update $\texttt{CurrentMax}(c)$ to $\texttt{NumSuff}(c)$.
\end{enumerate}
    
    \item After all iterations are finished, we have $\texttt{CurrentMax}(c) = \max_{d\in \cdisc{\cC}{A}} |\suff{d}{c}|$ for all $c \in \Ob\cC$. 
    For each $c \in \Ob\cC \setminus \Ob\cdisc{\cC}{A}$ we add  $c$ to $\rra{\cC}{A}$ if $2 \leq \texttt{CurrentMax}(c) $.
    The cost to check this is $O(n)$.
\end{enumerate}
    
    The total cost of step (2) across all $d\in \Ob\cdisc{\cC}{A}$ is 
    $O(|\Ob\cdisc{\cC}{A}|(n+m+|V|+|E|))$.
    Therefore, the total cost of the algorithm is
    \[ O(n+m |V|)+ O(|V||E|) + O((n+m+|V|+|E|)|\Ob\cdisc{\cC}{A}|)\]
    Since $|\Ob\cdisc{\cC}{A}| \leq |V| =p \leq n$, this is bounded above by $O((n+m)p)$.
\end{proof}

\begin{rmk}
    \label{rmk:gap}
There are important classes of rings for which projective presentations can be computed effectively.
For finite-dimensional algebras over a field, minimal projective resolutions of finitely generated modules can be computed by linear algebra, using an algorithm of Green--Solberg--Zacharia \cite{GSZ01}. 
For quotients of path algebras, Green and Solberg \cite{GS07} give an algorithmic construction of projective resolutions using noncommutative Gröbner-basis methods, although the resulting resolution need not be minimal in general. Both algorithms are implemented in the GAP package QPA \cite{QPA}.
These general-purpose algorithms do not exploit the additional structure of the constant functor $\intMod{\cC}$. This suggests the possibility for more specialized algorithms for computing a minimal initial subcategory.
\end{rmk}

\appendix
\section{Background Results From Representation Theory}\label{sec-background-rep-theory}

\begin{dfn}
A morphism $\psi:M\to N$ is \emph{right minimal} if every endomorphism $\mu: M \to M$ satisfying $\psi \mu=\psi$ is an automorphism. 
A free presentation $P_\bullet \to M$ as in \cref{eq:free-presentation} is \emph{minimal} if $\partial_0$ and $\partial_1$ are both right minimal. 
A \emph{projective cover} of $M$ is a right minimal epimorphism $\partial:P\to M$ such that $P$ is projective.
\end{dfn}
    
The following is a well known feature of minimal presentations.
\begin{lem}\label{lem:minimal-comparison}
Let $P_\bullet \to M$ be a minimal free presentation and let $Q_\bullet \to M$ be any free presentation. 
There exist retractions $r_0 : Q_0 \to P_0$ and $r_1 : Q_1 \to P_1$ such that the following diagram commutes
\[\begin{tikzcd}
Q_1 \arrow[d, "r_1"', dashed] \arrow[r, "q_1"] & Q_0 \arrow[d, "r_0"', dashed] \arrow[r, "q_0"] & M \arrow[d, "\id"'] \arrow[r] & 0 \\
P_1 \arrow[r, "p_1"]                           & P_0 \arrow[r, "p_0"]                           & M \arrow[r]                   & 0\,. 
\end{tikzcd}\]
\end{lem}

\begin{proof}
Projectivity (see \cref{prop:free-generators}) gives maps $s_0':P_0 \to Q_0$ and $r_0 : Q_0\to P_0$ such that $q_0s_0' = p_0$ and $p_0r_0 = q_0$. 
Hence, $p_0r_0s_0' = p_0$ and right minimality implies that $r_0s_0'$ is an automorphism. 
Then, $s_0 \defeq s_0'(r_0s_0')^{-1}$ is a section for $r_0$, so $r_0$ is the desired retraction.

Let $\overline{p_1}:P_1 \to \ker p_0$ and $\overline{q_1}:Q_1 \to \ker q_0$ denote the epimorphisms induced by taking images. 
Since $p_1$ is right minimal and $\ker p_0\hookrightarrow P_0$ is monic, the induced epimorphism $\overline p_1:P_1\twoheadrightarrow\ker p_0$ is right minimal. 
Both $s_0$ and $r_0$ restrict to maps $\overline{s_0} : \ker p_0 \to \ker q_0$ and $\overline{r_0} : \ker q_0 \to \ker p_0$. 
Since $P_1$ is projective and $\overline{q_1}$ is an epimorphism, there exists a lift $s_1' : P_1 \to Q_1$ such that 
$\overline{q_1}s_1'=\overline{s_0}\,\overline{p_1}$. 
Similarly, we get a map $r_1 : Q_1 \to P_1$ such that $\overline{p_1}r_1 = \overline{r_0}\,\overline{q_1}$. 
We compute $\overline{p_1}r_1s_1' = \overline{r_0}\,\overline{q_1}s_1' = \overline{r_0}\,\overline{s_0}\,\overline{p_1} = \overline{p_1}$ and right minimality makes $r_1s_1'$ invertible. 
We proceed as with $r_0$ to conclude the proof.
\end{proof}

The following allows one to apply certain results from module over rings to functor categories, and is the backbone of \cref{thm:proj-gamma}.
Item (\ref{item:GP-usual}) is an immediate consequence of the Gabriel-Popescu theorem \cite{GP64}.
Item (\ref{item:GP-equiv}) is Theorem IV.4.1 in \cite{Mitchell65}, but we give a proof based on the projectivization techniques in \cite[Chapter~II]{ARS}.

\begin{prp}\label{prp:projectivization}
    Let $\cC$ be a small category. Define
\begin{equation*}
    G\defeq \bigoplus_{c\in \Ob\cC} h_c \qquad \text{and}\qquad \Gamma \defeq \operatorname{End}_{\Mod\cC} (G)^\op \,.
\end{equation*}
Let $\Mod\Gamma$ denote the category of left modules over the ring $\Gamma$. Then
\begin{enumerate}
    \item\label{item:GP-usual} $\Hom(G,-):\Mod\cC\to \Mod \Gamma$ is a fully faithful functor and has an exact left adjoint $L:\Mod \Gamma \to \Mod \cC$.
    \item\label{item:GP-exactness} Since $G$ is projective $\Hom(G,-)$ is exact.
    \item\label{item:GP-weibel} $L$ preserves projective presentations.
    \item\label{item:GP-equiv} 
    If $\Ob\cC$ is finite, then 
    $\Hom(G,-)$ is an adjoint equivalence of categories.
\end{enumerate}
\end{prp}

\begin{proof} 
(\ref{item:GP-usual}) This is the usual statement of the Gabriel-Popescu theorem.
See, for instance, \cite[Theorem~8.5.8]{kashiwara05}.

(\ref{item:GP-exactness}) Since $G$ is projective, $\Hom(G,-)$ commutes with cokernels.

(\ref{item:GP-weibel}) The left adjoint $L$ is right exact and, since $\Hom(G,-)$ is exact, $L$ preserves projectives (see \cite[Prop.~2.3.10]{weibel94}).

(\ref{item:GP-equiv}) 
Write $H\defeq \Hom(G,-)$.
We first show $H$ is an equivalence of categories between projectives.
By \cref{prop:free-generators}, $G$ is a finite direct sum of finitely presentable functors, so $H$ commutes with colimits.
Suppose that $Q$ is projective.
Then for some set $A$, $Q$ is a direct summand of
\[X= \bigoplus_{a\in A}\Gamma = \bigoplus_{a\in A} HG = H \left(\bigoplus_{a\in A}G\right)\,.\]
Then $Q$ is the image of an idempotent $e:X \to X$, and since $H$ is full, $e=H\hat e$.
It follows that the image of $\hat e$ is projective and $Q \cong H(\im \hat e)$.
Thus, $H$ is essentially surjective onto the projective $\Gamma$-modules, and hence an equivalence, since $H$ is fully faithful.

It remains to show $H$ is essentially surjective.
Every $N \in \Ob\Mod\Gamma$ admits a projective presentation 
\[Q_1 \xrightarrow{q_1} Q_0 \xrightarrow{q_0} N\to 0\,.\]
By the above argument $H$ is an equivalence on projectives, so $q_1 \cong Hp_1$ for some morphism $p_1$. 
Then $N\cong \coker q_1 \cong H(\coker p_1)$, since $H$ is exact. 
\end{proof}

\cref{prp:projectivization} allows us to apply the following characterization of \emph{perfect rings}, due to Bass, to determine sufficient conditions for every $\cC$-module to have a projective cover.
This is the core principle behind \cref{cor:local-existence}.

For a module $M$ over a ring $\Gamma$, the \emph{Jacobson radical} of $M$, denoted $\rad M$, is the intersection of all \emph{maximal submodules} of $M$.

\begin{thm}[{\cite[Theorem~28.4]{AndersonFuller}}]\label{thm:bass}
    With the notation of \cref{prp:projectivization}, 
    the following are equivalent:
\begin{enumerate}[(i)]
    \item Every left $\Gamma$-module has a projective cover.
    \item $\Gamma/\rad \Gamma$ is semisimple and every non-zero left $\Gamma$-module contains a maximal submodule;
    \item $\Gamma$ contains no infinite orthogonal set of idempotents and every non-zero right $\Gamma$-module contains a minimal submodule.
\end{enumerate}
\end{thm}

The proof of \cref{thm:least-support} uses two classical decomposition theorems. 
The first one, \cref{lem:KRSA-free}, gives unique decomposition for free modules under the locality assumptions in \cref{thm:least-support}. 
It is a special case of Azumaya's generalization of the Krull--Remak--Schmidt theorem; see, for instance, \cite[Chapter~7,~Theorem~12.6]{AndersonFuller}. 
The second result, \cref{prp:projective-free}, is adapted from Warfield's version of a theorem of Crawley--J\'onsson; see \cite[Theorem~1]{warfield}. 
It shows that every projective $\cC$-module is isomorphic to a free module under the same locality assumption.

Both results were originally proved for modules over a ring. However, they have since been greatly generalized to sufficiently nice additive categories, including $\Mod\cC$; see \cite[Theorems~3 and~5]{WW76}.

\begin{lem}\label{lem:KRSA-free}
Suppose $\cC$ is a small category such that
$\kfield\cC(c,c)$ is a local ring for all $c\in\Ob\cC$.
Then every free $\cC$-module has a unique object support up to equivalence of categories.
\end{lem}

\begin{proof}
    Yoneda gives $\End_{\Mod\cC}(h_c) \cong \kfield\cC(c,c)^\op$,
    which is local by assumption.
    Let
\[
M=\bigoplus_{a\in A}h_{\sigma a}
\qquad\text{and}\qquad
M\cong\bigoplus_{b\in B}h_{\tau b}
\]
be two decompositions of a free module $M$ in $\Mod\cC$. 
By the Krull--Remak--Schmidt--Azumaya theorem,  
there is a bijection $\phi:A\to B$ such that
$h_{\sigma a}\cong h_{\tau(\phi a)}$
for every $a\in A$. 
By Yoneda, $\sigma a\cong\tau(\phi a)$, 
so the full subcategories on $\sigma A$ and $\tau B$ are equivalent. 
\end{proof}

\begin{prp}\label{prp:projective-free}
Suppose $\cC$ is a small category such that
$\kfield\cC(c,c)$ is a local ring for all $c\in\Ob\cC$.
If $P$ is a retract of a free module $M=\bigoplus_{a\in A}h_{\sigma a}$, then there exists a subset $A'\subseteq A$ giving a decomposition $P\cong\bigoplus_{a\in A'}h_{\sigma a}$.
\end{prp}

\begin{proof}
$P$ is isomorphic to a direct summand of $M$ by assumption.
By \cref{prop:free-generators}, each $h_{\sigma a}$ is finitely generated, and hence countably generated, and has local endomorphism ring. 
By \cite[Theorem~5]{WW76}, 
$P$ is a direct sum
of modules, each isomorphic to one of the $h_{\sigma a}$. 
\end{proof}

We conclude this appendix with the lemmas needed for \cref{sec:alt-thm-Ci}. 
Lemmas \ref{lem:quot-precomp-is-ff} and \ref{lem:quot-decomp-lift} are special cases of \cref{lem:lan-representable}.
\cref{lem:quot-pres-lift} is the strategy for the construction in \cref{sec:alt-thm-Ci}.

\begin{lem}\label{lem:quot-precomp-is-ff}
    Let $\pi:\cC\to \cD$ be a full and essentially surjective functor.
    Then
    $\pi^*: \Mod{\cD} \to \Mod\cC$ is fully faithful.
    Moreover, the unit $\eta:\id_{\Mod\cC} \Rightarrow \pi^*  \Lan_\pi$ is an epimorphism.
\end{lem}

\begin{proof}
    The fact that $\pi^*$ is fully faithful is well known; see for example \cite[Proposition~A.2]{vespa08}.
    Let $c,d\in \Ob\cC$.
    We can identify $\pi^*\Lan_\pi \kfield\cC(c, -)\cong\kfield\cD(\pi c, \pi -)$, by \cref{lem:lan-representable}.
    Then the component of the unit $(\eta_{h_c})_d:\kfield\cC(c,d) \twoheadrightarrow \kfield\cD(\pi c, \pi d)$ is given by the surjection $\gamma \mapsto \pi \gamma$, for all $\gamma\in \cC(c,d)$.
\end{proof}

\begin{lem}\label{lem:quot-decomp-lift}
Let $F:\cC\to \cD$ be an essentially surjective functor.
Suppose we have a function $\sigma:A \to \Ob\cD$ and a decomposition 
\[Q \cong\bigoplus_{a\in A}\kfield\cD(\sigma a,-)\]
of a free $\cD$-module $Q$.
Then there exists a function $\tau:A\to \Ob\cC$ and a decomposition 
\[ P =\bigoplus_{a\in A}\kfield\cC(\tau a,-)\]
of a free $\cC$-module $P$ such that $\sigma a \cong F \tau a$ for all $a\in A$
In particular, $Q \cong \Lan_F P$.
\end{lem}
\begin{proof}
    Since $F$ is essentially surjective, 
    we may choose for each $a\in A$ an object $\tau a \in \Ob\cC$ such that $\sigma a \cong F \tau a$. 
    This gives the desired function $\tau:A\to \Ob\cC$.
    By \cref{lem:lan-representable}, $\kfield\cD(F \tau a, -) \cong \Lan_F \kfield\cC(\tau a, -)$ for all $a\in A$.
    Since $\Lan_F$ commutes with coproducts, we have $Q \cong \Lan_F P$.
\end{proof}

\begin{lem}\label{lem:quot-pres-lift}
    Let $\pi:\cC\to \cD$ be a full and essentially surjective functor.
Let $P_0$ and $P_1$ be free $\cC$-modules such that
\[
\Lan_\pi P_1\xrightarrow{q_1}\Lan_\pi P_0\xrightarrow{q_0}N\longrightarrow0
\]
is a free presentation in $\Mod{\cD}$.
Let $\partial_1':P_1' \to P_0$ be a morphism such that $P_1'$ is free and its image is the kernel of the unit $\eta_{P_0}:P_0\twoheadrightarrow\pi^*\Lan_\pi P_0$.
Then there is a free presentation 
\[
P_1\oplus P_1'
\xrightarrow{[\partial_1~\partial_1']}
P_0
\xrightarrow{\partial_0} \pi^*N\longrightarrow0.
\]
\end{lem}

\begin{proof}
Since $\pi^*$ is exact, applying the functor to the free presentation $\Lan_\pi P_\bullet \to N$ yields an exact sequence. 
By \cref{lem:quot-precomp-is-ff}, the components of the unit
$
\eta_{P_n}:P_n\twoheadrightarrow\pi^*\Lan_\pi P_n
$
are epimorphisms.
Since $P_1$ is projective, there exists a lift $\partial_1:P_1 \to P_0$ of $(\pi^*q_1)\circ\eta_{P_1}$ such that the following commutes in $\Mod\cC$, where $\partial_0 \defeq (\pi^*q_0) \circ \eta_{P_0}$:
\begin{equation}
\begin{tikzcd}
P_1 \arrow[d, "\eta_{P_1}"', two heads] \arrow[r, "\exists \partial_1", dashed] & P_0 \arrow[d, "\eta_{P_0}"', two heads] \arrow[rd, "\partial_0"] &                  &   \\
\pi^*\Lan_\pi P_1 \arrow[r, "\pi^*q_1"']                                        & \pi^*\Lan_\pi P_0 \arrow[r, "\pi^*q_0"']              & \pi^*N \arrow[r] & 0.
\end{tikzcd}
\label{eq:dfn-partial1}
\end{equation}
Note that $\partial_0$ is clearly an epimorphism.

Since $\eta_{P_0}$ is an epimorphism, the snake lemma gives a short exact sequence (see \cite[Corollary~1.5]{triangle-lem} for details)
\begin{equation}\label{eq:sesKer}
0 \to \ker \eta_{P_0} \to \ker \partial_0 \xrightarrow{\eta_{P_0}} \ker(\pi^*q_0) \to 0.
\end{equation}
Since $\eta_{P_1}$ is an epimorphism commutativity of \cref{eq:dfn-partial1} gives
\[\im(\eta_{P_0}\partial_1)=
\im((\pi^*q_1)\circ\eta_{P_1})=\im(\pi^*q_1)=\ker(\pi^*q_0)\,.\]
In particular the induced morphism $\overline{\partial_1}:P_1 \to \ker(\pi^*q_0)$ is an epimorphism.
Now choose an epimorphism $\partial_1':P_1' \to \ker \eta_{P_0}$ with $P_1'$ a free module.
Applying the \emph{horseshoe lemma} (see \cite[Lemma~2.2.8]{weibel94}) to \cref{eq:sesKer} yields an epimorphism 
\[[\partial_1'~ \partial_1]:P_1'\oplus P_1 \to\ker \partial_0
,\]
which gives the desired presentation.
\end{proof}

\bibliographystyle{halpha}
\bibliography{functors}

\end{document}